\documentstyle[11pt]{article}
\begin{document}
\renewcommand{\thefootnote}{\fnsymbol{footnote}}
\pagestyle{empty}
\setcounter{page}{1}
\newfont{\twelvemsb}{msbm10 scaled\magstep1}
\newfont{\eightmsb}{msbm8} \newfont{\sixmsb}{msbm6} \newfam\msbfam
\textfont\msbfam=\twelvemsb \scriptfont\msbfam=\eightmsb
\scriptscriptfont\msbfam=\sixmsb \catcode`\@=11
\def\Bbb{\ifmmode\let\next\Bbb@\else \def\next{\errmessage{Use
      \string\Bbb\space only in math mode}}\fi\next}
\def\Bbb@#1{{\Bbb@@{#1}}} \def\Bbb@@#1{\fam\msbfam#1}
\newfont{\twelvegoth}{eufm10 scaled\magstep1}
\newfont{\tengoth}{eufm10} \newfont{\eightgoth}{eufm8}
\newfont{\sixgoth}{eufm6} \newfam\gothfam
\textfont\gothfam=\twelvegoth \scriptfont\gothfam=\eightgoth
\scriptscriptfont\gothfam=\sixgoth \def\frak{\frak@}
\def\frak@#1{{\fam\gothfam{{#1}}}} \def\frak@@#1{\fam\gothfam#1}
\catcode`@=12
%
%
%
\def\CC{{\Bbb C}}
\def\NN{{\Bbb N}}
\def\QQ{{\Bbb Q}}
\def\RR{{\Bbb R}}
\def\ZZ{{\Bbb Z}}
\def\cA{{\cal A}}          \def\cB{{\cal B}}          \def\cC{{\cal C}}
\def\cD{{\cal D}}          \def\cE{{\cal E}}          \def\cF{{\cal F}}
\def\cG{{\cal G}}          \def\cH{{\cal H}}          \def\cI{{\cal I}}
\def\cJ{{\cal J}}          \def\cK{{\cal K}}          \def\cL{{\cal L}} 
\def\cM{{\cal M}}          \def\cN{{\cal N}}          \def\cO{{\cal O}}
\def\cP{{\cal P}}          \def\cQ{{\cal Q}}          \def\cR{{\cal R}} 
\def\cS{{\cal S}}          \def\cT{{\cal T}}          \def\cU{{\cal U}}
\def\cV{{\cal V}}          \def\cW{{\cal W}}          \def\cX{{\cal X}}
\def\cY{{\cal Y}}          \def\cZ{{\cal Z}}
\def\qed{\hfill \rule{5pt}{5pt}}
\def\arcsinh{\mathop{\rm arcsinh}\nolimits}
\newtheorem{theorem}{Theorem}
\newtheorem{prop}{Proposition}
\newtheorem{conj}{Conjecture}
\newenvironment{result}{\vspace{.2cm} \em}{\vspace{.2cm}}

\pagestyle{plain}
\begin{center}
$\phantom{xx}$
\smallskip
\smallskip

{\LARGE {\bf JORDANIAN QUANTUM $\left(\hbox{SUPER}\right)$ALGEBRAS ${\cal U}_{\sf h}\left(
{\sf g}\right)$ VIA CONTRACTION METHOD AND MAPPING: REVIEW}} \\[0.5cm]

{B. ABDESSELAM$^{a}$, R. CHAKRABARTI$^{b}$, A. YANALLAH$^{c}$ and M.B ZAHAF$^{d}$}

\smallskip 

$^{a,c,d}${\it Laboratoire de Physique 
Quantique de la Mati\`ere et Mod\'elisations Math\'ematiques (LPQ3M), Centre Universitaire de 
Mascara, 29000-Mascara, Alg\'erie}

{\it $^{a}$Laboratoire de Physique Th\'eorique d'Oran, Université d'Oran Es-S\'enia, 
31100-Oran, Alg\'erie}

{\it $^{b}$Department of Theoretical Physics, University of Madras, Guindy Campus, Madras 
600025, India}
 
\end{center}


\begin{abstract}
{\small Recenty, a class of transformations of ${\cal R}_{q}$-matrices was introduced such 
that the $q\longrightarrow 1$ limit gives explicit nonstandard ${\cal R}_{h}$-matrices. The 
transformation matrix is singular as $q\longrightarrow 1$. For the transformed matrix, the 
singularities, however, cancel yielding a well-defind construction. We have shown that our 
method can be implemented systematically on ${\cal R}_{q}$-matrices of all dimensions of 
${\cal U}_q(sl(N))$, ${\cal U}_q(osp(1|2))$ and ${\cal U}_q(sl(2|1))$. Explicit 
constructions are presented for ${\cal U}_q(sl(2))$, ${\cal U}_q(sl(3))$, ${\cal U}_q
(osp(1|2))$ and ${\cal U}_q(sl(2|1))$, while choosing ${\cal R}_{q}$ for (fund. rep.) 
$\otimes$(arbitray irrep.). Our method yields nonstandard deformations along with a nonlinear 
map of the ${\sf h}$-Borel subalgebra on the corresponding classical Borel subalgebra, which 
can be easily extended to the whole algebra. Following this approach, we construct explicitly 
here the jordanian quantum (super)algebras (nonstandard version) ${\cal U}_{h}(sl(2))$, 
${\cal U}_{h}(sl(3))$, ${\cal U}_{h}(osp(1|2))$ and ${\cal U}_{h}(sl(2|1))$. These 
Hopf (super)algebras are equipped with a remarkably simpler coalgebraic 
structure. Generalizing 
our results on ${\cal U}_{h}(sl(3))$, we give the higher dimensional Jordanian 
(super)algebras ${\cal U}_{h}(sl(N))$ for all $N$. The universal ${\cal R}_{h}$-matrices 
are also given.} 
\end{abstract}

\tableofcontents

\section{Introduction and motivation}

It is well known that any enveloping (super)algebra ${\cal U}({\sf g})$ of an Lie 
(super)algebra ${\sf g}$ has many quantizations: the first one called the {\it Drinfeld-Jimbo 
deformation or the standard quantum deformation} \cite{P1,P2} is quasitriangular, whereas the 
others one called the {\it Jordanian deformations or the nonstandard quantum deformations} 
\cite{P3} are triangular $\left({\cal R}_{21}{\cal R}=1\right)$. A typical example of the 
Jordanian quantum algebras was first introduced by Ohn \cite{P4}. In general, nonstandard 
quantum algebras are obtained by applying twist \cite{P5} to the corresponding Lie algebras. 
A twisting that produces an algebra isomorphism to the Ohn algebra \cite{P4} is found in 
\cite{P6,P7}. 

Recently, the twisting procedure was extensively employed to study a wide variety of Jordanian 
deformed algebras, such as ${\cal U}_{h}(sl(N))$ algebra [8-12], symplectic algebra 
${\cal U}_{h}(sp(N))$ \cite{P13}, orthogonal algebra ${\cal U}_{h}(so(N))$ [14-17] and 
orthosymplectic superalgebra ${\cal U}_{h}(osp(1|2))$ [18-20]. It follows from these 
studies that:

\begin{itemize}

\item The nonstandard quantum algebras have undeformed commutation relations;

\item The Jordanian deformation appears only in the coalgebraic structure; 

\item  the coproduct and the antipode maps have very complicated forms in comparison with the 
Drinfeld-Jimbo and the Ohn deformations. 

\end{itemize}

So far Jordanian quantum algebra ${\cal U}_{h}(sl(N))$ has been explicitly written, with a 
simple coalgebra but with deformed commutation relations, only for $N=2$ \cite{P4}. This 
amounts to a choice of an appropriate basis, in which the commutation relations are deformed
but the corresponding coalgebraic structure remains simple. The main object here  
is to construct somes nonstandard versions of an enveloping Lie algebra which have deformed 
commutation relations; but are endowed with a relatively simpler coalgebraic structure compared 
to those in the previous studies [8-20].      

Following this approach, we have proposed, recently [21-25], a new technic which makes possible the 
construction of a nonstandard version ${\cal U}_{\sf h}({\sf g})$ of an enveloping Lie 
(super)algebra ${\cal U}({\sf g})$ by a suitable contraction, from the standard ones ${\cal U}_q
({\sf g})$. Our scheme consist to obtain the ${\cal R}_{\sf h}$-matrix, for all 
dimensions of a (super)Jordanian quantum (super)algebra ${\cal U}_{\sf h}({\sf g})$ from 
the ${\cal R}_q$-matrix associated to the standard quantum (super)algebra ${\cal U}_q
({\sf g})$ through a specific transformation ${\sf G}$ (singular in the $q\longrightarrow 1$), 
as follows:
\begin{equation}
{\cal R}_{\sf h}=\lim_{q\longrightarrow 1}\left[{\sf G}^{-1}\otimes {\sf G}^{-1}\right]
{\cal R}_q\left[{\sf G}\otimes {\sf G}\right],
\end{equation}          
where, for example, ${\sf G}={\sf E}_q\left(\frac{{\sf h}E_{1N}}{q-1}\right)$ for ${\cal U}_q
\left(sl\left(N\right)\right)$ ($E_{1N}$ is the longer positive root generator of the Lie 
algebra $sl(N)$), ${\sf G}={\sf E}_q\left(\frac{{\sf h}e^2}{q^2-1}\right)$ for ${\cal U}_q
\left(osp\left(1|2\right)\right)$ ($e$ is the fermionic positive simple root of the Lie algebra
$osp(2|1)$) and ${\sf G}={\sf E}_q\left(\frac{{\sf h}e_{1}}{q-1}\right)$ for ${\cal U}
\left(sl\left(2|1\right)\right)$ ($e_1$ is the bosonic positive simple root generator of the 
Lie algebra $sl(2|1)$). The deformed exponential map ${\sf E}_q$ is defined by
\begin{equation}
{\sf E}_q\left(\eta\right)=\sum_{n=0}^\infty\frac{\left(\eta\right)^n}{\left[n\right]_q!},
\qquad [n]_q=\frac{q^{n}-q^{-n}}{q-q^{-1}},\qquad [n]_q!=[n]_q\times [n-1]_q!,\qquad [0]_q!=1.
\end{equation}    
This procedure yields nonstandard deformation along with a nonlinear map of the 
${\sf h}$-Borel subalgebra on the corresponding classical Borel subalgebra, which can be easily 
extended to the whole algebra. Following here this strategy, we have construct the Jordanian 
quantum (super)algebras ${\cal U}_{\sf h}(sl(2))$, ${\cal U}_{\sf h}(sl(3))$, ${\cal U}_{\sf h}
(osp(1|2))$ and ${\cal U}_{\sf h}(sl(2|1))$, wherein we use the contraction procedure and the 
map mentionned above. The nonstandard versions obtained here have deformed 
commutation relations, but the coalgebraic part is more simple and in compact form.  

The manuscrit is organized as follows: In section 2, the jordanian quantum algebra 
${\cal U}_{\sf h}\left(sl\left(2\right)\right)$ and ${\cal U}_{\sf h}\left(sl\left(3\right)
\right)$ are constructed. A nonlinear map between 
${\cal U}_{\sf h}\left(sl\left(2\right)\right)$ and ${\cal U}\left(sl\left(2\right)\right)$ 
(resp. ${\cal U}_{\sf h}\left(sl\left(3\right)\right)$ and ${\cal U}\left(sl\left(3\right)
\right)$) is then established. Generalizing our results on the ${\cal U}\left(sl\left(3\right)
\right)$ algebra, the higher dimensional algebras ${\cal U}_{\sf h}\left(sl\left(N\right)
\right)$, $N>4$, are also introduced via a nonlinear map and proved to be a Hopf algebra 
endowed with a triangular ${\cal R}_{\sf h}$-matrix. The (super)jordanian quantum 
versions ${\cal U}_{\sf h}\left(sl(2|1)\right)$ and ${\cal U}_{\sf h}\left(osp\left(1|2\right)
\right)$ are presented respectively in section 3 and 4. Finally, we conclude in section 5.     
             
\section{${\cal U}_{{\sf h}}(sl(N))$: Map, Hopf Algebra, Irreps. and 
${\cal R}_{{\sf h}}$-matrix}
 
For our purpose, the deformation parameter ${\sf h}$ is an arbitrary complex number. For 
simplicity, we start with ${\cal U}_q(sl(2))$ algebra.    

\subsection{The ${\cal R}_{\sf h}$-matrices of the Nonstandard ${\cal U}_{\sf h}(sl(2))$ 
Algebra}

The generating elements $\left( H,\,X,\,Y\,\right) $ of the algebra ${\cal U}_{\sf h}\left(sl
\left(2\right)\right)$ obey the following commutation rules \cite{P4} 
\begin{eqnarray}
&&\left[H,X\right]=2\frac{\sinh {\sf h}X}{\sf h},\qquad
\left[H,Y\right]=-Y\left(\cosh {\sf h}X\right)-\left(\cosh {\sf h}X\right)Y,\qquad
\left[X,Y\right]=H.
\end{eqnarray}
The non-commutative coproduct structure of ${\cal U}_{\sf h}(sl(2))$ read \cite{P4} 
\begin{eqnarray}
&&\Delta_{\sf h}\left(X\right)=X\otimes 1+1\otimes X,\qquad
\Delta_{\sf h}\left(Y\right)=Y\otimes e^{{\sf h}X}+e^{-{\sf h}X}\otimes Y,\qquad
\Delta_{\sf h}\left(H\right)=H\otimes e^{{\sf h}X}+e^{-{\sf h}X}\otimes H. 
\end{eqnarray}
The Borel subalgebra generated by the elements $\left(H,\,X\right)$ is equiped with a Hopf 
structure. The universal ${\cal R}$-matrix for ${\cal U}_{\sf h}\left(sl\left(2\right)\right)$ 
may be cast in the form \cite{P26} 
\begin{equation}
{\cal R}_{\sf h}=\exp\biggl(-{\sf h}X\otimes e^{{\sf h}X}H\biggr)\exp\biggl({\sf h}e^{{\sf h}
X}H\otimes X\biggr)
\end{equation}
which obviously coincides with the universal ${\cal R}$-matrix of the borel subalgebra. The 
fundamental $(j=\frac{1}{2})$ representation of the algebra (3) remains undeformed \cite{P4}: 
\begin{equation}
X=\left(
\begin{array}{cc} 0 & 1 \\ 0 & 0
\end{array}
\right),\,\,\,\,\,\,\,Y=\left( 
\begin{array}{cc}
0 & 0 \\ 
1 & 0
\end{array}
\right),\,\,\,\,\,\,\,H=\left( 
\begin{array}{cc}
1 & 0 \\ 
0 & -1
\end{array}
\right). 
\end{equation}
Using the fundamental representation (6) in the first sector of the tensor product of the 
operators in the expression (5) of the universal ${\cal R}_{\sf h}$-matrix, we obtain the 
$R_{\sf h}$-matrix in the $\left(\frac{1}{2}\otimes j\right)$ representation: 
\begin{equation}
R_{\sf h}=\left(\begin{array}{cc}
e^{{\sf h}X} & -{\sf h}H+\frac{\sf h}{2}\left(e^{{\sf h}X}-e^{-{\sf h}X}\right) \\ 
0 & e^{-{\sf h}X}
\end{array}\right).
\end{equation}
Absent from upper triangular form (7) and indeed from the universal ${\cal R}_{\sf h}$-matrix 
(5) is the generator $Y$ completing the ${\cal U}_{\sf h}\left(sl\left(2\right)\right)$ algebra.
We will show how (7) can be obtained, directly and for arbitrary $j$, from the corresponding 
$R_{q}$-matrix for $\left(\frac{1}{2}\otimes j\right)$ representation given by (see, for 
example \cite{P27}) 
\begin{equation}
{\em R}_{q}=\left( 
\begin{array}{cc}
q^{{\cal J}_0/2} & q^{\frac{1}{2}}\left(1-q^{-2}\right) {\cal J}_- \\ 
0 & q^{-{\cal J}_0/2}
\end{array}\right).
\end{equation}
Here the generators of ${\cal U}_{q}\left(sl\left(2\right)\right)$ are denoted by $\left(q^{\pm 
{\cal J}_{0}},{\cal J}_{\pm}\right)$ satisfying the standard relations 
\begin{equation}
q^{{\cal J}_{0}}{\cal J}_{\pm}={\cal J}_{\pm }q^{{\cal J}_{0}\pm 2},\,\,\,\,\,\,\,\,\,\left[
{\cal J}_{+},\,{\cal J}_{-}\right] =\frac{q^{{\cal J}_{0}}-q^{-{\cal J}_{0}}}{q-q^{-1}}
\equiv \left[{\cal J}_{0}\right]. 
\end{equation}
Hereafter we consider generic $q$, excluding roots of unity.

For the purpose of transformating $R_{q}$-matrix in (8), we now consider a $q$-deformed 
exponential operator: ${\sf E}_{q}\left(\eta {\cal J}_{+}\right)=\sum_{n=0}^{\infty}\frac{\left(\eta 
{\cal J}_{+}\right)^n}{\left[n\right]!}$. We choose, for an arbitrary finite constant ${\sf h}$ 
the parameter $\eta$ as $\eta=\frac{\sf h}{q-1}$. We emphasize that through the deformed 
exponential ${\sf E}_q(x)$ defined in (2) has non convenient simple expression for its inverse 
(comparable to the standard $q$-exponential \cite{P27} satisfying $\exp_q(x)^{-1}=\exp_{q^{-1}}(-x)$), 
this is precisely what is needed for obtaining non-singular limiting forms for the $R$-matrix 
elements for arbitrary representations and other interesting properties. For any given value of 
$j$, the series (2) may be terminated after setting ${\cal J}_+^{2j+1}=0$; but, we proceed 
quite generally as follows. Defining 
\begin{equation}
{\cal T}_{(\alpha)}= {\sf E}^{-1}_q\left(\eta {\cal J}_+\right){\sf E}_q\left(q^\alpha\eta 
{\cal J}_+\right) 
\end{equation}
with ${\cal T}_{(0)}=1$, we obtain 
\begin{equation}
{\sf E}^{-1}_q\left(\eta{\cal J}_+\right)q^{\alpha{\cal J}_0/2}{\sf E}_q\left(q^\alpha\eta 
{\cal J}_+\right)={\cal T}_{(\alpha)}q^{\alpha{\cal J}_0/2}.
\end{equation}
For transforming $R_q$-matrix in (8), the operators ${\cal T}_{\pm 1}$ are of particular 
importance. We will be concerned with simple rational values of $\alpha$. For later use, we 
note the identity ${\sf E}^{-1}_q\left(\eta{\cal J}_+\right)q^{(\alpha+\beta){\cal J}_0/2}
{\sf E}_q\left(q^\alpha\eta{\cal J}_+\right)=\left({\sf E}^{-1}_q\left(\eta{\cal J}_+\right)
q^{\alpha{\cal J}_0/2}{\sf E}_q\left(q^\alpha\eta{\cal J}_+\right)\right)\left({\sf E}^{-1}_q
\left(\eta{\cal J}_+\right)q^{\beta{\cal J}_0/2}{\sf E}_q\left(q^\alpha\eta{\cal J}_+\right)
\right)$, which in notation (11) reads 
\begin{equation}
{\cal T}_{(\alpha+\beta)}q^{(\alpha+\beta){\cal J}_0/2}=\left({\cal T}_{(\alpha)}
q^{\alpha{\cal J}_0/2}\right)\left({\cal T}_{(\alpha)}q^{\alpha{\cal J}_0/2}\right). 
\end{equation}  
Moreover using the identiy $\left[{\cal J}_+^n,{\cal J}_-\right]=\frac{[n]}
{q-1}\left(q^{{\cal J}_0/2}{\cal J}_+^{n-1}q^{{\cal J}_0/2}-q^{-{\cal J}_0/2}{\cal J}_+^{n-1}
q^{-{\cal J}_0/2}\right)$, the following commutator is obtained $\left[{\sf E}_q\left(\eta{\cal J}_+
\right),{\cal J}_-\right]=\frac{[n]}{q-1}\left(E_q\left(q\eta{\cal J}_+\right)q^{{\cal J}_0}-
E_q\left(q^{-1}\eta{\cal J}_+\right)q^{-{\cal J}_0}\right)$, which, in turn, leads to 
$E_q\left(\eta{\cal J}_+\right)^{-1}{\cal J}_-$ $E_q\left(\eta{\cal J}_+\right)=-\frac{\eta}{q-
q^{-1}}\left({\cal T}_{(1)}q^{{\cal J}_0}-{\cal T}_{(-1)}q^{-{\cal J}_0}\right)+
{\cal J}_-$. Evaluating term by term, the $q\longrightarrow 1$ limits of ${\cal T}_{(\pm 1)}$
are found to be finite and of the form
\begin{equation}
\lim_{q\rightarrow 1}{\cal T}_{(\pm 1)}=T_{(\pm 1)}=\sum _{n=0}^\infty c_n^{(\pm)}
\left({\sf h}J_+\right)^n,
\end{equation}
where $\left(J_0,J_{\pm}\right)$ are the generators of the classical $sl(2)$ algebra ($\left[J_0,
J_{\pm}\right]=\pm 2J_{\pm},\;[J_+,J_-]=J_0$). The first few coefficients $\{c_n^{(\pm)}|n\geq 
0\}$ in (13) read $c_0^{(\pm)}=1$, $c_1^{(\pm)}=\pm 1$, $c_2^{(\pm)}=1/2$, $c_3^{(\pm)}=0$, 
$c_4^{(\pm)}=-1/8$, $c_5^{(\pm)}=0$. If the limits are indeed finite, then from (12) it is 
evident that 
\begin{equation}
T_{(\alpha)}=\left(T_{(1)}\right)^\alpha,
\end{equation}
where $T_{(\alpha)}=\lim_{q\rightarrow 1}{\cal T}_{(\alpha)}$. Henceforth we write
$T_{(1)}=T$. 
The result obtained here suggests the following derivation of the closed form of $T$. To this 
end, we left and right multiply the second commutaion relation in (9) by ${\sf E}_q^{-1}\left(
\eta{\cal J}_+\right)$ and ${\sf E}_q\left(\eta{\cal J}_+\right)$, respectively ${\sf E}_q^{-1}\left(\eta{\cal J}_+
\right)\left(q^{{\cal J}_0}-q^{-{\cal J}_0}\right){\sf E}_q\left(\eta{\cal J}_+\right)=
\left(q-q^{-1}\right){\sf E}_q^{-1}\left(\eta{\cal J}_+\right)\left[{\cal J}_+, {\cal J}_-
\right]$ ${\sf E}_q\left(\eta{\cal J}_+\right)$. Using (10), (11) and (12), we obtain 
\begin{equation}
{\cal T}_{(2)}q^{{\cal J}_0}-{\cal T}_{(-2)}q^{-{\cal J}_0}={\sf h}(q+1)\left(
{\cal T}_{(1)}{\cal J}_+q^{{\cal J}_0}+q^{-2}{\cal T}_{(-1)}{\cal J}_+q^{-{\cal J}_0}\right)+
\left( q^{{\cal J}_0}-q^{-{\cal J}_0}\right). 
\end{equation}
Using (14), we now obtain the following equation for $T$ as $q\longrightarrow 1$: $T^2-T^{-2}
=\left(T+T^{-1}\right)\left(2{\sf h}J_+\right)$, which after a factorization yields 
\begin{equation}
T-T^{-1}=2{\sf h}J_+.
\end{equation}
The quadratic relation (16) in $T$ is now solved
  \begin{equation}
T^{\pm 1}=\pm {\sf h}J_++\sqrt{ 1+{\sf h}^2J_+^2}.
\end{equation}
This our crucial result. 

With all these results now in hand we go back to $R_q$ in (8). We choose the transformation 
matrix as $G=g_{\frac{1}{2}}\otimes g$, where $g={\sf E}_{q}\left( \eta J_{+}\right)$ and 
$g_{\frac{1}{2}}\equiv g\left|_{j=\frac{1}{2}}\right.$. We obtain 
\begin{eqnarray}
R_{\sf h}&=&\lim_{q\rightarrow 1}\left(G^{-1}R_{q}G\right)=\left( 
\begin{array}{cc}
T & -\frac{\sf h}{2}\left( T+T^{-1}\right)J_{0}+\frac{\sf h}{2}\left(T-T^{-1}\right) \\ 
0 & T^{-1}
\end{array}
\right),  
\end{eqnarray}
where we have defined 
\begin{eqnarray}
&&e^{\pm {\sf h}X} =T^{\pm 1}=\pm {\sf h}J_++\left(1+{\sf h}^2J_{+}^2
\right)^{\frac{1}{2}},\qquad
H=\frac{1}{2}\left(T+T^{-1}\right)J_{0}=\left(1+{\sf h}^2J^2_{+}\right)J_0. 
\end{eqnarray}
It is easy to verify that 
\begin{equation}
\left[H,T^{\pm 1}\right]=T^{\pm 2}-1\Rightarrow \left[H,X\right]=2\frac{\sinh {\sf h}X}{\sf h}. 
\end{equation}
Comparing (7) with (18), we see that the contraction scheme, which comprises our transformation 
and limiting procedure has furnished the $R_{\sf h}$-matrix along with a nonlinear map of the 
subalgebra of ${\cal U}_{\sf h}(sl(2))$ generated by $\left(H,X\right)$ on the classical one 
generated by ($J_0,\,J_+$). Indeed, defining 
\begin{equation}
Y=J_-\frac{{\sf h}^2}{4}J_+\left(J_0^{2}-1\right)
\end{equation}
we complete the ${\cal U}_{\sf h}(sl(2))$ algebra and show that 
\begin{equation}
\left[H,Y\right]=-Y\left(\cosh {\sf h}X\right)-\left(\cosh{\sf h}X\right)Y,\qquad 
\left[X,Y\right]=H.
\end{equation}
The expressions (19) and (21) may be looked as a particular realization of the ${\cal U}_{{\sf h}}
(sl(2))$ generators $(H,X,Y)$. We have therefore, developed an invertible nonlinear map of $(H,X,Y)$
on classical generators $(J_0,J_+,J_-)$.  

\subsection{The ${\cal R}_{\sf h}$-matrices of the Nonstandard ${\cal U}_{\sf h}(sl(3))$ Algebra}

The major interest of our method is that can be generalized for obtaining the nonstandard 
$R_{\sf h}$-matrices of algebras of higher dimensions as contraction limits of the 
corresponding $R_q$-matrices. Here we treat the $sl(3)$ algebra. Choosing the Chevalley 
generators corresponding to the simple roots of the ${\cal U}_{q}(sl(3))$ algebra as 
$k_1^{\pm}=q^{\pm h_1}$, $k_2^{\pm}=q^{\pm h_2}$, $k_3^{\pm}=q^{\pm h_3}=q^{\pm (h_1+h_2)}$, 
${\hat e}_1$, ${\hat e}_2$, ${\hat e}_3={\hat e}_1 {\hat e}_2-q^{-1}{\hat e}_2{\hat e}_1$, 
${\hat f}_1$, ${\hat f}_2$ and ${\hat f}_3={\hat f}_2{\hat f}_1-q{\hat f}_1{\hat f}_2$. The 
Hopf structure of the ${\cal U}_q(sl(3))$ algebra is defined by [27]
\begin{eqnarray}
&& [h_i,\;h_j]=0,\qquad  [h_i,{\hat e}_j]=a_{ij}{\hat e}_j,\qquad 
[h_i,{\hat f}_j]=-a_{ij}{\hat f}_j,\qquad [{\hat e}_i,{\hat f}_j]=\delta_{ij}[h_i],\nonumber\\
&& {\hat e}_1{\hat e}_{3}=q{\hat e}_3{\hat e}_1,\qquad {\hat e}_2{\hat e}_{3}=q^{-1}{\hat e}_3
{\hat e}_2,\qquad {\hat f}_1{\hat f}_{3}=q{\hat f}_3{\hat f}_1,\qquad {\hat f}_2{\hat f}_{3}=
q^{-1}{\hat f}_3{\hat f}_2,\\
&&\Delta_{q}\left(h_i\right) =h_i\otimes 1+1\otimes h_i,\qquad 
\Delta_{q}\left({\hat e}_i\right) ={\hat e}_i\otimes q^{h_i/2}+q^{-h_i/2}\otimes {\hat e}_i,
\qquad 
\Delta_{q}\left({\hat f}_i\right) ={\hat f}_i\otimes q^{h_i/2}+q^{-h_i/2}\otimes {\hat f}_i,
\nonumber
\end{eqnarray}          
where $[\chi]=\frac{q^\chi-q^{-\chi}}{q-q^{-1}}$. The Cartan matrix for the $sl(3)$ algebra 
reads $a=\pmatrix{2& -1\cr \\ -1& 2}$. The universal matrix of the ${\cal U}_{\sf h}(sl(3))$ 
algebra is given by
\begin{eqnarray}
 {\cal R}_q&=&q^{\sum_{i,j}(a^{-1})_{ij}h_i\otimes h_j}\exp_{q^{-2}}\left(\lambda {\hat e}_2
q^{h_2/2}\otimes q^{-h_2/2}{\hat f}_2\right)\exp_{q^{-2}}\left(\lambda {\hat e}_3
q^{h_3/2}\otimes q^{-h_3/2}{\hat f}_3\right)\nonumber\\
&&\exp_{q^{-2}}\left(\lambda {\hat e}_1 q^{h_1/2}\otimes q^{-h_1/2}{\hat f}_1\right),
\end{eqnarray}
where $\lambda=q-q^{-1}$, $\exp_q\left(x\right)=\sum_{n=0}^\infty x^n/\{n\}_q!$, $\{n\}_q=
\{n\}_q\{n-1\}_q!$, $\{0\}_q=1$ and $\{n\}_q=\frac{1-q^n}{1-q}$. We denote the classcial 
generators ($q=1$) of the $sl(3)$ algebra by $h_1$, $h_2$, $h_3=h_1+h_2$, $e_1$, 
$e_2$, $e_3=e_1e_2-e_2e_1$, $f_1$, $f_2$ and $f_3=f_2f_1-f_1f_2$.

Let us turn now to the nonstandard $R_{\sf h}$-matrices. For brevity and simplicity 
we limit here our self to (fundamental irrep.)$\;\otimes\;$(arbitrary irrep.). Recall that for 
${\cal U}_q(sl(3))$ algebra the $R_q$-matrix in the representation (fund.)$\;\otimes\;$(arb.) 
reads \cite{P27}: 
\begin{eqnarray}
&& R_q=\biggl(\pi_{\hbox{fund.}}\;\otimes\;\pi_{\hbox{arb.}}\biggr) {\cal R}_q
= \pmatrix{q^{{1\over 3}(2h_1+h_2)} & q^{{1\over 3}(2h_1+h_2)}\Lambda_{12} & q^{{1\over 3}
(2h_1+h_2)}\Lambda_{13} \cr 0 & q^{-{1\over 3}(h_1-h_2)} & q^{-{1\over 3}(h_1-h_2)}
\Lambda_{23} \cr 0 & 0 &q^{-{1\over 3}(h_1+2h_2)} \cr},
\end{eqnarray}
where 
\begin{eqnarray}
&&\Lambda_{12}=q^{-1/2}\bigl(q-q^{-1}\bigr)q^{-h_1/2}{\hat f}_1,\,\,
\Lambda_{13}=q^{-1/2}\bigl(q-q^{-1}\bigr){\hat f}_3 q^{-{\frac{1}{2}}(h_1+h_2)},\,\, 
\Lambda_{23}=q^{-1/2}\bigl(q-q^{-1}\bigr)q^{-h_2/2}{\hat f}_2.
\end{eqnarray}

We have shown in \cite{P21} that the nonstandard $R_{\sf h}$-matrix (in the representation 
(fund.)$\;\otimes\;$(arb.)) arise from the $R_q$-matrix as follows: 
\begin{eqnarray}
&& R_{\sf h}=\lim_{q\rightarrow 1} \biggl[{\sf E}_q\biggl({\frac{{\sf h}{\hat e}_3}{q-1}}
\biggr)_{\hbox{(fund.)}}\otimes {\sf E}_q\biggl({\frac{{\sf h}{\hat e}_3}{q-1}}
\biggr)_{\hbox{(arb.)}}\biggr]^{-1} R_q\biggl[{\sf E}_q\biggl({\frac{{\sf h}{\hat e}_3}{q-1}}
\biggr)_{\hbox{(fund.)}} \otimes {\sf E}_q\biggl({\frac{{\sf h}{\hat e}_3}{q-1}}
\biggr)_{\hbox{(arb.)}}\biggr]\nonumber\\
&&\phantom{R_{\sf h}}=\lim_{q\rightarrow 1}\pmatrix{ {\sf E}_q^{-1}\bigl({{\sf h}
{\hat e}_3\over q-1}\bigr) & 0 & -{{\sf h}{\sf E}_q^{-1}\bigl({{\sf h}{\hat e}_3
\over q-1}\bigr)\over q-1} \cr 0 & {\sf E}_q^{-1}\bigl({{\sf h}{\hat e}_3\over q-1}\bigr)& 0\cr 
0 & 0 &{\sf E}_q^{-1}\bigl({{\sf h}{\hat e}_3\over q-1}\bigr)\cr} R_q \pmatrix{
{\sf E}_q\bigl({{\sf h}{\hat e}_3\over q-1}\bigr) & 0 & {{\sf h}{\sf E}_q\bigl({{\sf h}
{\hat e}_3\over q-1}\bigr)\over q-1} \cr 0 & {\sf E}_q\bigl({{\sf h}{\hat e}_3\over q-1}\bigr)& 
0\cr 0 & 0 & {\sf E}_q\bigl({{\sf h}{\hat e}_3\over q-1}\bigr)\cr}  \nonumber \\
&& \phantom{R_{\sf h}}=\pmatrix{T &2{\sf h}T^{-1/2}e_2 & -{{\sf h}\over
2}(T+T^{-1})\bigl(h_1+h_2\bigr) +{{\sf h}\over 2}\bigl(T-T^{-1}\bigr)\cr 0 &
I & -2{\sf h}T^{1/2}e_1 \cr 0 & 0 & T^{-1} \cr},
\end{eqnarray}
where 
\begin{eqnarray}
&& T={\sf h}e_3+\sqrt{1+{\sf h}^2e_3^2}, \qquad\qquad T^{-1}=-{\sf h}e_3+
\sqrt{1+{\sf h}^2e_3^2}.
\end{eqnarray}

The following properties are pointed out:

\smallskip

{\bf 1.} The corner elements of (27) have exactly the same structure as in 
the $R_{\sf h}$-matrix of ${\cal U}_{\sf h}(sl(2))$. This implies
that the classical generators $e_3$, $h_3=h_1+h_2$ and $f_3$ of ${\cal U}(sl(3))$ 
are deformed (for the nonstandard quantization: ${\cal U}(sl(3))\longrightarrow 
{\cal U}_{\sf h}(sl(3))$) as follows \cite{P21}: 
\begin{eqnarray}
&& T^{\pm 1}=\pm {\sf h}e_3+\sqrt{1+{\sf h}^2e_3^2}, \qquad
H_3=\sqrt{1+{\sf h}^2e_3^2}h_3, \qquad F_3=f_3-\frac{{\sf h}^2}{4}e_3
\bigl(h_3^2-1\bigr),
\end{eqnarray}
and satisfy evidently the commutation relations \cite{P4} 
\begin{eqnarray}
&& TT^{-1}=T^{-1}T=1,\qquad \qquad [H_3,T^{\pm 1}]=T^{\pm 2}-1, \nonumber\\
&&[T^{\pm 1},F_3]=\pm {\frac{{\sf {h}}}{2}}\biggl(H_3T^{\pm 1}+T^{\pm 1}H_3\biggr),\qquad 
[H_3,F_3]=-{\frac{1}{2}}\biggl(TF_3+F_3T+T^{-1}F_3 +F_3T^{-1}\biggr).
\end{eqnarray}
With the following definition (see Ref. \cite{P4}) $E_3={\sf h}^{-1}\ln T={\sf h}^{-1}
\hbox{arcsinh}\;{\sf h}e_3$, it follows that the elements $H_3$, $E_3$ and $F_3$ satisfy 
the relations of the ${\cal U}_{\sf h}(sl(2))$ algebra \cite{P4}
\begin{eqnarray}
&& [H_3,E_3]=2\frac{\sinh {\sf h}E_3}{\sf h},\qquad [H_3,F_3]=-F_3
\biggl(\cosh {\sf h}E_3\biggr)- \biggl(\cosh {\sf h}E_3\biggr)F_3,\qquad
[E_3,F_3]=H_3,
\end{eqnarray}
where it is obvious that as ${\sf h}\longrightarrow 0$ , we have $(H_3,E_3,F_3)
\longrightarrow (h_3,e_3,f_3)$. It is now evident from (31) that 
${\cal U}_{\sf h}(sl(2))\subset {\cal U}_{\sf h}(sl(3))$.

\smallskip

{\bf 2.} The expression (27) of the $R_{\sf h}$-matrix indicate that the simple root generators 
$e_1$ and $e_2$ are deformed as follows: 
\begin{eqnarray}
&& E_1=\sqrt{{\sf h}e_3+\sqrt{1+{\sf h}^2e_3^2}}e_1=T^{1/2}e_1,  \qquad 
E_2=\sqrt{{\sf h}e_3+\sqrt{1+{\sf h}^2e_3^2}}e_2=T^{1/2}e_2.
\end{eqnarray}
To complete our algebra ${\cal U}_{{\sf h}}(sl(3))$, let us introduce the
following ${\sf h}$-deformed generators: 
\begin{eqnarray}
&& F_1=\sqrt{-{\sf h}e_3+\sqrt{1+{\sf h}^2e_3^2}}f_1 +{\frac{{\sf h}}{2}}%
\sqrt{{\sf h}e_3 +\sqrt{1+{\sf h}^2e_3^2}}e_2h_3= T^{-1/2}\biggl(f_1+{\frac{%
{\sf h}}{2}}e_2Th_3\biggr),  \nonumber \\
&& F_2=\sqrt{-{\sf h}e_3+\sqrt{1+{\sf h}^2e_3^2}}f_2 -{\frac{{\sf h}}{2}}%
\sqrt{{\sf h}e_3 +\sqrt{1+{\sf h}^2e_3^2}}e_1h_3= T^{-1/2}\biggl(f_2-{\frac{%
{\sf h}}{2}}e_1Th_3\biggr),  \nonumber \\
&& H_1=\biggl(-{\sf h}e_3+\sqrt{1+{\sf h}^2e_3^2}\biggr) \biggl(\sqrt{1+{\sf %
h}^2e_3^2}h_1 +{\frac{{\sf h}}{2}}e_3(h_1-h_2)\biggr)= h_1-{\frac{{\sf h}}{2}%
}e_3T^{-1}h_3,  \nonumber \\
&& H_2=\biggl(-{\sf h}e_3+\sqrt{1+{\sf h}^2e_3^2}\biggr) \biggl(\sqrt{1+{\sf %
h}^2e_3^2}h_2 -{\frac{{\sf h}}{2}}e_3(h_1-h_2)\biggr)= h_2-{\frac{{\sf h}}{2}%
}e_3T^{-1}h_3.
\end{eqnarray}
The expressions (29), (32) and (33) constitute a realization of the Jordanian
algebra ${\cal U}_{\sf h}(sl(3))$ with the classical generators via a
nonlinear map. This immediately yields the irreducible representations
(irreps.) of ${\cal U}_{\sf h}(sl(3))$ in an explicit and simple manner.

\begin{prop}
The Jordanian algebra ${\cal U}_{\sf h}(sl(3))$ is then an associative algebra over $\CC$ 
generated by $H_1$, $H_2$, $H_3$, $E_1$, $E_2$, $T$, $T^{-1}$, $F_1$, $F_2$ 
and $F_3$, satisfying the commutation relations [22]
\begin{eqnarray}
&& [H_1,H_2]=0, \qquad\qquad
     [H_1,T^{-1}H_3]=[H_2,T^{-1}H_3]=0,\nonumber \\
&& [H_1,E_1] =2E_1,\qquad  [H_2,E_2]=2E_2,\qquad [H_1,E_2]=-E_2,
\qquad [H_2,E_1]=-E_1, \nonumber\\ 
&& [T^{-1}H_3,E_1]=E_1,\qquad [T^{-1}H_3,E_2]=E_2, \qquad [H_1,F_1]=-2F_1+{\sf h}E_2T^{-1}H_3,
\nonumber \\
&& [H_2,F_2] =-2F_2-{\sf h}E_1T^{-1}H_3, \qquad 
[H_1,F_2] =F_2-{\sf h}E_1T^{-1}H_3, \qquad [H_2,F_1] =F_1+{\sf h}E_2T^{-1}H_3, \nonumber\\
&& [TH_3,F_1] =-T^2F_1, \qquad\qquad\qquad\qquad\qquad\;\; 
     [TH_3,F_2] =-T^2F_2, \nonumber\\
&& [T^{-1}E_1,F_1]={1\over 2}(T+T^{-1})H_1+{1\over 2}(T-T^{-1})H_2, \qquad
[T^{-1}E_2,F_2]={1\over 2}(T+T^{-1})H_2+{1\over 2}(T-T^{-1})H_1, \nonumber\\ 
&& [T^{-1}E_1,F_2] =0,\qquad [T^{-1}E_2,F_1] =0, \qquad 
[E_1,E_2] ={1\over 2{\sf h}}(T^2-1),\nonumber\\
  && [TF_2,TF_1]=T\biggl(F_3-{{\sf h}\over 2}H_3TH_3
                 -{{\sf h}\over 8}(T-T^{-1})\biggr) \nonumber\\
  && [TH_1,T^{\pm 1}]= {1\over 2}(T^{\pm 2}-1),\qquad  
 [TH_2,T^{\pm 1}]= {1\over 2}(T^{\pm 2}-1),\nonumber\\
  && [H_1,F_3]= -{T^{-1}\over 4}\biggl(TF_3+F_3T+T^{-1}F_3
                +F_3T^{-1}\biggl)-{{\sf h}\over 4}T^{-1}H_3^2
                -{{\sf h} \over 4}H_3T^{-1}H_3, \nonumber \\
  && [H_2,F_3]= -{T^{-1}\over 4}\biggl(TF_3+F_3T+T^{-1}F_3
                +F_3T^{-1}\biggl)-{{\sf h}\over 4}T^{-1}H_3^2
                -{{\sf h} \over 4}H_3T^{-1}H_3, \nonumber \\
  && [E_1,T]=[E_1,T^{-1}]=[E_2,T]=[E_2,T^{-1}]=0,\nonumber\\
  && [F_1,T]= {\sf h}TE_2, \qquad [F_1,T^{-1}]= -{\sf h}T^{-1}E_2, \qquad 
[F_2,T]= -{\sf h}TE_1, \qquad [F_2,T^{-1}]={\sf h}T^{-1}E_1,\nonumber\\
  && [E_1,F_3]=-{1\over 2}\biggl(TF_2+F_2T\biggr), \qquad\qquad\;\;\;\;\;
     [E_2,F_3]={1\over 2}\biggl(TF_1+F_1T\biggr), \nonumber \\
&& [F_1,F_3]={\sf h}TF_1-{\sf h}E_2F_3+{{\sf h}^2\over 4}TE_2,\qquad 
   [F_2,F_3]={\sf h}TF_2+{\sf h}E_1F_3-{{\sf h}^2\over 4}TE_1.
\end{eqnarray}
\end{prop}

We have stated here the final results only. To obtain the expressions of $H_1$ and $H_2$, 
we have proceed as follows: By analogy with (29), we have first started with the ansatz 
$\sqrt{1+{\sf {h}}^2e_3^2}h_1$ and $\sqrt{1+{\sf {h}}^2e_3^2}h_2$. It is easy to see that 
$[\sqrt{1+{\sf {h}}^2e_3^2}h_1,F_3]=-{\frac{1}{4}}\biggl(TF_3+F_3T+T^{-1}F_3 +F_3T^{-1}\biggr)
+{\frac{{\sf {h}}^2}{4}}\biggl(e_3(h_1-h_2)H_3 +H_3e_3(h_1-h_2)\biggr)$ and 
$[\sqrt{1+{\sf {h}}^2e_3^2}h_2,F_3]=-{\frac{1}{4}}\biggl(TF_3+F_3T+T^{-1}F_3 +F_3T^{-1}\biggr)
-{\frac{{\sf {h}}^2}{4}}\biggl(e_3(h_1-h_2)H_3+H_3e_3(h_1-h_2)\biggr)$.
Then, if we add to $\sqrt{1+{\sf {h}}^2e_3^2}h_1$ and we deduct from $\sqrt{1+{\sf {h}}^2
e_3^2}h_2$ the term ${\frac{{\sf {h}}}{2}}e_3(h_1-h_2)$, we obtain 
$[\sqrt{1+{\sf {h}}^2e_3^2}h_1+{\frac{{\sf {h}}}{2}}e_3(h_1-h_2),F_3] =-{\frac{1}{4}}
\biggl(TF_3+F_3T+T^{-1}F_3+F_3T^{-1}\biggr)
+{\frac{{\sf {h}}}{4}}T(h_1-h_2)H_3+{\frac{{\sf {h}}}{4}}H_3T(h_1-h_2)$ 
and , 
$[\sqrt{1+{\sf {h}}^2e_3^2}h_2-{\frac{{\sf {h}}}{2}}e_3(h_1-h_2),F_3] =-{
\frac{1}{4}}\biggl(TF_3+F_3T+T^{-1}F_3+F_3T^{-1}\biggr)
-{\frac{{\sf {h}}}{4}}T(h_1-h_2)H_3-{\frac{{\sf {h}}}{4}}H_3T(h_1-h_2)$.
These commutation relations suggest to take $H_1\sim\biggl(\sqrt{1+{\sf {h}}
^2 e_3^2}h_1+{\frac{{\sf {h}} }{2}}e_3(h_1-h_2)\biggr)$ and $H_2\sim\biggl(%
\sqrt{1+ {\sf {h}}^2e_3^2}h_2-{\frac{{\sf {h}} }{2}}e_3(h_1-h_2)\biggr)$.
Finally, to preserve the Cartan subalgebra, we are obliged to multiply $%
\biggl(\sqrt{1+{\sf {h}}^2e_3^2}h_1 +{\frac{{\sf {h}} }{2}}e_3(h_1-h_2)%
\biggr)$ and $\biggl(\sqrt{1+{\sf {h}}^2e_3^2} h_2-{\frac{{\sf {h}} }{2}}%
e_3(h_1-h_2)\biggr)$ respectively by $T^{-1}$, i.e. to take $H_1=T^{-1}%
\biggl(\sqrt{1+{\sf {h}}^2e_3^2}h_1+{\frac{{\sf {h}} }{2}}e_3 (h_1-h_2)%
\biggr)=h_1-{\frac{{\sf h}}{2}}e_3T^{-1}h_3$ and $H_2=T^{-1}\biggl
(\sqrt{1+{\sf {h}}^2e_3^2}h_2-{\frac{{\sf {h}} }{2}}e_3(h_1-h_2)\biggr)=h_2- 
{\frac{{\sf h}}{2}}e_3T^{-1}h_3$. The expressions of $F_1$ and $F_2$ are
obtained in similair way. The expressions (29), (32) and (33) may be looked
now as a particular realization of the ${\cal U}_{{\sf h}}(sl(3))$
generators. Others maps can be also considered.

In terms of the Chevalley generators (simple roots) $\{E_1, E_2,F_1, 
F_2,H_1,H_2\}$, the algebra ${\cal U}_{{\sf{h}}}(sl(3))$ is defined as follows [22]:
\begin{eqnarray}
&&T=\biggl(1+2{\sf h}[E_1,E_2]\biggr)^{1/2},\qquad\qquad\qquad\;\;
T^{-1}=\biggl(1+2{\sf h}[E_1,E_2]\biggr)^{-1/2},\nonumber\\
&&[H_1,H_2]=0, \qquad [H_1,E_1] =2E_1,\qquad [H_2,E_2] =2E_2,\qquad [H_1,E_2]=-E_2,
\nonumber\\
&& [H_2,E_1]=-E_1, \qquad [H_1,F_1] =-2F_1+{\sf h}E_2(H_1+H_2), \qquad
[H_2,F_2] =-2F_2-{\sf h}E_1(H_1+H_2), \nonumber\\
&& [H_1,F_2] =F_2-{\sf h}E_1(H_1+H_2), \qquad\qquad\;
     [H_2,F_1] =F_1+{\sf h}E_2(H_1+H_2), \nonumber\\
&& [T^{-1}E_1,F_1]={1\over 2}(T+T^{-1})H_1+{1\over 2}(T-T^{-1})H_2,\qquad 
   [T^{-1}E_2,F_2]={1\over 2}(T+T^{-1})H_2+{1\over 2}(T-T^{-1})H_1,\nonumber\\
&& [T^{-1}E_1,F_2] =[T^{-1}E_2,F_1] =0, \nonumber \\
&& E_1^2E_2-2E_1E_2E_1+E_2E_1^2=0,\qquad
   E_2^2E_1-2E_2E_1E_2+E_1E_2^2=0,\nonumber \\
&& (TF_1)^2TF_2-2TF_1TF_2TF_1+TF_2(TF_1)^2=0,
   \qquad
   (TF_2)^2TF_1-2TF_2TF_1TF_2+TF_1(TF_2)^2=0,
\end{eqnarray}
or, briefly
\begin{eqnarray}
&& [H_i,H_j]=0, \qquad
[H_i,E_j] =a_{ij}E_j,\qquad [H_i,F_j] =-a_{ij}F_j+T^{-1}[F_j,T](H_1+H_2),\nonumber \\
&& [T^{-1}E_i,F_j]=\delta_{ij}\biggl(T^{-1}H_i+{1\over 2}
(T-T^{-1})(H_1+H_2)\biggr), \nonumber\\ 
&& (\hbox{ad}\; E_i)^{1-a_{ij}}(E_j)=(\hbox{ad}\; TF_i)^{1-a_{ij}}(TF_j)=0,\qquad i\neq j,
\end{eqnarray}
where $(a_{ij})_{i,j=1,2}$ is the Cartan matrix of $sl(3)$. Let us turn now to the coalgebra 
structure:

\begin{prop}
The Jordanian quantum 
algebra ${\cal U}_{\sf h}(sl(3))$ admits a Hopf structure with coproducts, 
antipodes and counits determined by [22]
\begin{eqnarray}
&&\Delta (E_1)=E_1\otimes 1+T\otimes E_1,\qquad 
\Delta (E_2)=E_2\otimes 1+T\otimes E_2, \qquad 
\Delta (T^{\pm 1})=T^{\pm 1}\otimes T^{\pm 1},\nonumber \\
&&\Delta(F_1)=F_1\otimes 1+T^{-1}\otimes F_1+{\sf h}H_3\otimes E_2,\qquad
\Delta(F_2)=F_2\otimes 1+T^{-1}\otimes F_2-{\sf{h}}H_3\otimes E_1,\nonumber\\
&&\Delta (F_3)=F_3\otimes T+T^{-1}\otimes F_3, \qquad
\Delta (H_1)=H_1\otimes 1+1\otimes H_1-{1\over 2}(1-T^{-2})\otimes T^{-1}H_3, \nonumber \\
&&\Delta (H_2)=H_2\otimes 1+1\otimes H_2-{1\over 2}(1-T^{-2})\otimes T^{-1}H_3,\qquad
\Delta (H_3)=H_3\otimes T+T^{-1}\otimes H_3,  \nonumber \\
&& S(E_1)=-T^{-1}E_1,\qquad  S(E_2)=-T^{-1}E_2,\qquad
   S(T)=T^{-1},\qquad S(T^{-1})=T,\nonumber \\
&& S(F_1)=-TF_1+{\sf h}TH_3T^{-1}E_2,\qquad S(F_2)=-TF_2-{\sf{h}}TH_3T^{-1}E_1,\nonumber \\
&& S(F_3)=-TF_3T^{-1},\qquad S(H_1)=-H_1-{1\over 2}(T-T^{-1})H_3, \nonumber \\
&& S(H_2)=-H_2-{1\over 2}(T-T^{-1})H_3, \qquad S(H_3)=-TH_3T^{-1},\nonumber\\
&&\epsilon(a)=0,\qquad\forall 
a\in\biggl\{H_1,H_2,H_3,E_1,E_2,F_1,F_2,F_3\biggr \},\nonumber \\
&& \epsilon (T)=\epsilon (T^{-1})=1.  
\end{eqnarray}
\end{prop}

All the Hopf algebra axioms can be verified by direct calculations. Let us remark that our 
coproducts have simpler forms compared to those maps in [8-20]. This is one benifit of our 
procedure. Pertinent to the algebraic structures of our Hopf algebra described
in (30), (34) and (37), here we obtain its universal ${\cal R}$-matrix in the following form
[22]: 
\begin{eqnarray}
&& {\cal R}_{\sf h}=\exp\biggl(-{\sf h}E_3\otimes TH_3\biggr)\exp\biggl(2{\sf h}TH_3
\otimes E_3\biggr). 
\end{eqnarray}
The above univesal matrix satisfies the required properties for the full Hopf structure 
discussed earlier. We note that the element, generated by $E_3$ and $H_3$ coincides with the 
universal ${\cal R}_{\sf h}$-matrix of the subalgebra involving the generators corresponding 
to the highest root, and may be connected to the results obtained by the contraction process 
by a suitable twist operator that can be derived as a series expansion in ${\sf h}$. 
The nonstandard algebras ${\cal U}_{\sf h}(sl(4))$ and ${\cal U}_{\sf h}(sl(5))$ can be 
derived in a similar way (see refs. \cite{P22}).  

\subsection{${\cal U}_{{\sf h}}(sl(N))$: Generalization}

>From these above studies, It is easy to see that:

\begin{prop}
The Jordanian quantization deform ${\cal U}(sl(N))$'s Chevalley generators as follows [22]: 
\begin{eqnarray}
&& T={\sf h}[e_1,[e_2,\cdots,[e_{N-2},e_{N-1}]]]+
\sqrt{1+{\sf h}^2([e_1,[e_2,\cdots,[e_{N-2},e_{N-1}]]])^2}, \nonumber \\
&& T^{-1}=-{\sf h}[e_1,[e_2,\cdots,[e_{N-2},e_{N-1}]]]+
\sqrt{1+{\sf h}^2([e_1,[e_2,\cdots,[e_{N-2},e_{N-1}]]])^2}, \nonumber \\ 
&&E_i=T^{(\delta_{i1}+\delta_{i,N-1})/2}e_i,\qquad i=1,\cdots,N-1,\nonumber\\
&& F_i=T^{-(\delta_{i1}+\delta_{i,N-1})/2}\biggl(
f_i+{{\sf h}\over 2}T[f_i,[e_1,[e_2,\cdots,[e_{N-2},e_{N-1}]]]](h_1+\cdots +
h_{N-1})\biggr)\nonumber\\
&& H_i=h_i-{(\delta_{i1}+\delta_{i,N-1}){\sf h}\over 2}
[e_1,[e_2,\cdots,[e_{N-2},e_{N-1}]]]T^{-1}(h_1+\cdots +
h_{N-1}) 
\end{eqnarray}
and they satisfy the commutation relations
\begin{eqnarray}
&& [H_i,H_j]=0, \qquad [H_i,E_j] =a_{ij}E_j,\nonumber \\
&& [H_i,F_j] =-a_{ij}F_j+(\delta_{i1}+\delta_{i,N-1})T^{-1}[F_j,T]
(H_1+\cdots + H_{N-1}), \nonumber\\
&& [T^{-(\delta_{i1}+\delta_{i,N-1})}E_i,F_j]=
\delta_{ij}\biggl(T^{-(\delta_{i1}+\delta_{i,N-1})}H_i
+{(\delta_{i1}+\delta_{i,N-1})
\over 2}(T-T^{-1})(H_1+\cdots +H_{N-1})\biggr),\nonumber\\
&& [E_i,E_j] =0,\qquad\qquad \qquad\qquad |i-j|>1,\nonumber \\ 
&& [T^{(\delta_{i1}+\delta_{i,N-1})}F_i,T^{(\delta_{j1}+\delta_{j,N-1})}F_j]=0,
\qquad\qquad |i-j|>1,\nonumber \\ 
&& (\hbox{ad} \;E_i)^{1-a_{ij}}(E_j)=
(\hbox{ad} \; T^{(\delta_{i1}+\delta_{i,N-1})}F_i)^{1-a_{ij}}
(T^{(\delta_{j1}+\delta_{j,N-1})}F_j)=0,\qquad (i\neq j),
\end{eqnarray} 
where $(a_{ij})_{i,j=1,\cdots , N}$ is the Cartan matrix of $sl(N)$, 
i.e. $a_{ii}=2$, $a_{i,i\pm 1}=-1$ and $a_{ij}=0$ for $|i-j|>1$.    
\end{prop}

The algebra (40) is called the Jordanian quantum algebra ${\cal U}_{\sf h}(sl(N))$. The 
expressions (39) may be regarded as a particular nonlinear realisation of the ${\cal U}_{\sf h}
(sl(N))$ generators. The Jordanian algebra ${\cal U}_{\sf h}(sl(N))$ (40) admits the following 
coalgebra structure [22]:  
\begin{eqnarray}
&&\Delta (E_i)=E_i\otimes 1+T^{(\delta_{i1}+\delta_{i,N-1})}
\otimes E_i, \nonumber \\
&&\Delta (F_i)=F_i\otimes 1+T^{-(\delta_{i1}+\delta_{i,N-1})} 
\otimes F_i+T(H_1+\cdots +H_{N-1})\otimes T^{-1}[F_i,T],\nonumber \\
&&\Delta (H_i)=H_i\otimes 1+1\otimes H_i-
{(\delta_{i1}+\delta_{i,N-1})\over 2}(1-T^{-2})\otimes (H_1+\cdots +H_{N-1}), 
\nonumber\\
&& S(E_i)=-T^{-(\delta_{i1}+\delta_{i,N-1})}E_i, \qquad
S(F_i)=-T^{(\delta_{i1}+\delta_{i,N-1})}F_i+T^2(H_1+\cdots +H_{N-1})
T^{-2}[F_i,T], \nonumber \\
&& S(H_i)=-H_i+{(\delta_{i1}+\delta_{i,N-1})\over 2}(1-T^2)
(H_1+\cdots +H_{N-1}), \nonumber \\
&& \epsilon(E_i)= \epsilon(F_i)=\epsilon(H_i)=0.
\end{eqnarray}

Following (38), we obtain the universal ${\cal R}_{\sf h}$-matrix of an arbitrary  
${\cal U}_{\sf h}(sl(N))$ algebra in the following general form [22]:
\begin{eqnarray}
&& {\cal R}_{\sf h}=\exp\biggl(-{\sf h}E_{1N}\otimes TH_{1N}\biggr)\exp\biggl(2{\sf h}TH_{1N}
\otimes E_{1N}\biggr). 
\end{eqnarray}
where $H_{1N}=T(H_1+\cdots +H_{N-1})$ and $E_{1N}={\sf h}\ln T$. The above universal 
${\cal R}_{\sf h}$-matrix of the full algebra Hopf algebra is obtained from the generators 
associated with the highest root; and its coincides with the univesal ${\cal R}_{\sf h}$-matrix 
of the Hopf subalgebra cite associated with the highest root. It is intersting to note that the 
nonlinear map (40) equips the ${\sf h}$-deformed generators $(E_{i},F_{i},H_{i})$ with an 
additional induced co-commutative coproduct. Simlilarly, the may be also equipped with 
an induced co-commutative coproduct. Similarly, the undeformed generators $(e_{i},f_{i},h_{i})$, 
via the inverse map, may be viewed as elements of the ${\cal U}_{\sf h}(sl(N))$ algebra; and, 
thus, may be endowed with an induced noncommutative coproduct.

\section{Nonstandard quantizations of $osp(1|2)$ superalgebra}

It has been recently demonstrated \cite{JS98} that three distinct bialgebra structures exist 
on the classical $osp(1|2)$ superalgebra, and all of them are coboundary. The classical Lie
superalgebra $osp(1|2)$ has three even $(h, b_{\pm})$ and two odd $(e, f)$ generators, which 
obey the commutation relations 
\begin{eqnarray}
&&[h,\;e] = e,\qquad [h,\;f] = -f,\qquad \{e,\;f\}=-h,  \qquad
[h,\;b_{\pm}] = {\pm}2 b_{\pm},\qquad \;[b_+,\;b_-] = h,  \nonumber \\
&&[b_+,\;f]=e, \qquad [b_-,\;e]=f,\qquad b_+ = e^2,\qquad b_- = - f^2.
\label{eq:clalg}
\end{eqnarray}
The generators $(h, b_{\pm})$ form a subalgebra $sl(2) \subset osp(1|2)$.
The inequivalent classical $r$-matrices defined on $osp(1|2)$ superalgebra
have been listed in Ref. \cite{JS98} as 
\begin{equation}
r_1 = h \wedge b_+,\qquad r_2 = h \wedge b_+ - e \wedge e,\qquad r_3 = t (h
\wedge b_+ + h \wedge b_- - e \wedge e - f \wedge f).  \label{eq:clrma}
\end{equation}
The standard quasi-triangular $q$-deformation of the $osp(1|2)$ superalgebra
considered in Refs. \cite{KR89, S90} corresponds to $r_3$. The parameter $t$
in $r_3$ becomes irrelevant in quantization as it can be absorbed into the
deformation parameter. The $r_1$ matrix is comprised of the elements of the $sl(2)$
subalgebra. This allows quantization of the $osp(1|2)$ superalgebra using
the inclusion $sl(2) \subset osp(1|2)$. This is done \cite{CK98} by applying
Drinfeld twist for the $sl(2)$ subalgebra to the full $osp(1|2)$
superalgebra. The Hopf algebra ${\bf U_h}(osp(1|2))$ obtained thereby is
triangular. The quantization of the $r_1$ matrix has been obtained in Ref. 
\cite{CK98} in terms of the classical basis set with undeformed commutation
relations but with coproduct structures deformed in a complicated manner.

Recently the classical matrix $r_2$ has been quantized \cite{ACS03} using
nonlinear basis elements. The corresponding quantized algebra ${\cal U}_{\sf h}
(osp(1|2))$ is known \cite{ACS03, BLT03} to satisfy the triangularity condition. An 
important issue observed \cite{K98} in this context is that the quantum $R_{\sf h}$ 
matrix of the ${\cal U}_{\sf h}(osp(1|2))$ algebra in the fundamental representation may be
obtained {\it via} a contraction mechanism from the corresponding $R_q$
matrix of the standard ${\cal U}_{q}(osp(1|2))$ algebra in the $q\rightarrow 1$ limit. A 
generalization of this contraction procedure for arbitrary representations, though clearly 
desirable as it will allow us to systematically obtain various quantities of interest of the 
${\cal U}_{\sf h}(osp(1|2))$ algebra from the corresponding quantities of the $q$
-deformed ${\cal U}_{q}(osp(1|2))$ algebra, has not been achieved so far.

In another problem considered also here, we express the other nonstandard
Hopf algebra ${\bf U_h}(osp(1|2))$ corresponding to the classical $r_1$
matrix in a {\it nonlinear} basis. Here we follow the approach in Ref. 
[21], where the Jordanian ${\cal U}_{\sf h}(sl(2))$ algebra
has been introduced in terms of a nonlinear basis set, while retaining the
coproduct structure of these basis elements simple. A consequence of our
choice of nonlinear basis elements is that, Ohn's ${\cal U}_{\sf h}(sl(2))$ algebra 
\cite{P4} explicitly arises as a Hopf subalgebra of our ${\bf U_h}(osp(1|2))$ 
algebra. This feature is not directly evident in the construction given in Ref. \cite{CK98}. 
Moreover, while the algebraic relations of our ${\bf U_h}(osp(1|2))$ algebra are 
deformed, {\it the coproduct structures are considerably simple}. Our approach may be of
consequence in building physical models of many-body systems employing
coalgebra symmetry \cite{BH99}. Invertible nonlinear maps, and the twist
operators pertaining to these maps, exist connecting the deformed and the
undeformed basis sets. We will present here a {\it class} of invertible
maps interrelating the Hopf algebra ${\bf U_h}(osp(1|2))$, based on
quantization of the $r_1$ matrix, with its classical analog ${\cal U}(osp(1|2))$. The twist 
operators \textit{vis-{\`a}-vis} the above maps will be discussed in the sequel. It is shown 
that a particular map called `minimal twist map' implements the simplest twist given directly 
by the factorized form of the universal ${\bf R_h}$ matrix of the ${\bf U_h}(osp(1|2))$ 
algebra. For a `non-minimal' map the twist has an additional factor. We evaluate this twist 
operator as a series in the deformation parameter.

Using the oft-used nomenclature we will refer to the algebra obtained by
quantizing the $r_2$ matrix as the super-Jordanian ${\cal U}_{\sf h}(osp(1|2))$ algebra, 
whereas the algebra generated \textit{via} quantization of $r_1$ matrix will be noted as 
Jordanian ${\bf U_h}(osp(1|2))$ algebra.

\subsection{Super-Jordanian ${\cal U}_{\sf h}(osp(1|2))$ algebra 
\textit{via} contraction process}

The Hopf structure of the super-Jordanian ${\cal U}_{\sf h}(osp(1|2))$
algebra using nonlinear basis elements was obtained \cite{ACS03}
previously. For comparing with our subsequent results we list it, after a
slightly altered normalization, as 
\begin{eqnarray}
&&[H, E] = \frac{1}{2}\,(T + T^{-1})\,E , \qquad [H, F] = - \frac{1}{4}\,(T
+ T^{-1})\,F - \frac{1}{4}\,F\,(T + T^{-1}),  \nonumber \\
&&\{E, F\} = - H,\qquad [H, T^{\pm 1}] = T^{\pm 2} - 1,  \nonumber \\
&&[H, Y] = - \frac{1}{2}(T + T^{-1}) Y - \frac{1}{2} Y (T + T^{-1}) - \frac{%
\mathsf{h}}{4} E (T - T^{-1}) F - \frac{\mathsf{h}}{4} F (T - T^{-1}) E, 
\nonumber \\
&&[T^{\pm 1},Y]=\pm\frac{\sf h}{2}(T^{\pm 1} H+HT^{\pm 1}),
\qquad E^2=\frac{T-T^{-1}}{2{\sf h}}, \qquad F^2 =-Y,  \nonumber \\
&&[T^{\pm 1}, F] = \pm {\sf h} T^{\pm 1} E,\qquad [Y, E] = \frac{1}{4} (T
+ T^{-1}) F + \frac{1}{4} F (T + T^{-1}),  \nonumber \\
&&\Delta(H) = H \otimes T^{-1} + T \otimes H + {\sf h} E T^{1/2} \otimes
E T^{-1/2},\qquad \Delta(E) = E \otimes T^{-1/2} + T^{1/2} \otimes E, 
\nonumber \\
&&\Delta(F) = F \otimes T^{-1/2} + T^{1/2} \otimes F, \qquad \Delta(T^{\pm
1}) = T^{\pm 1} \otimes T^{\pm 1},  \nonumber \\
&&\Delta(Y) = Y \otimes T^{-1} + T \otimes Y + \frac{\mathsf{h}}{2} E
T^{1/2} \otimes T^{-1/2} F + \frac{\mathsf{h}}{2} T^{1/2} F \otimes E
T^{-1/2},  \nonumber \\
&&\varepsilon(H) = \varepsilon(E) = \varepsilon(F) =\varepsilon(Y) = 0,
\qquad \varepsilon(T^{\pm 1}) = 1,  \nonumber \\
&&S(H) = - H - {\sf h} E^2, \qquad S(E) = - E, \qquad S(F) = - F + 
\frac{\sf h}{2} E,  \nonumber \\
&&S(T^{\pm 1}) = T^{\mp 1}, \qquad S(Y) = -Y +\frac{\sf h}{2} H + \frac{\sf h^2}{4} E^2,  
\label{eq:r2hopf}
\end{eqnarray}
where ${\sf h}$ is the deformation parameter. The ${\cal U}_{\sf h}(osp(1|2))$
algebra has only \textit{one} Borel subalgebra generated by the
elements $(H, E, T^{\pm 1})$. Kulish observed \cite{K98} that the $R_{
\mathsf{h}}$ matrix in the fundamental representation of the super-Jordanian 
${\cal U}_{\sf h}(osp(1|2))$ algebra may be obtained \textit{via} a
transformation, singular in the $q \to 1$ limit, from the corresponding $R_q$
matrix in the fundamental representation of the standard $q$-deformed 
${\cal U}_q(osp(1|2))$ algebra.

Our task here is to generalize the above contraction
procedure for arbitrary representations. As an application of our
contraction scheme, we construct the $L$ operator corresponding to the Borel
subalgebra of the super-Jordanian ${\cal U}_{\sf h}(osp(1|2))$
algebra from the corresponding $L$ operator of the standard $q$-deformed 
${\cal U}_{\sf h}(osp(1|2))$ algebra. To this end, we first quote some
well-known \cite{KR89, S90} results on the ${\cal U}_{\sf h}(osp(1|2))$
algebra. The ${\cal U}_q(osp(1|2))$ algebra is generated by three
elements $({\hat h},\;{\hat e},\;{\hat f})$ obeying the Hopf structure 
\begin{eqnarray}
&&[{\hat h},\;{\hat e}] = {\hat e},\qquad\qquad [{\hat h},\;{\hat f}] = - {%
\hat f},\qquad\qquad \{{\hat e},\;{\hat f}\}=- [h]_{q},  \nonumber \\
&&\Delta\left({\hat h}\right)={\hat h}\otimes 1+1\otimes {\hat h},\qquad
\Delta\left({\hat e}\right)={\hat e}\otimes q^{-{\hat h}/2} +q^{{\hat h}%
/2}\otimes {\hat e},\qquad \Delta\left({\hat f}\right)={\hat f}\otimes q^{-{%
\hat h}/2} +q^{{\hat h}/2}\otimes {\hat f},  \nonumber \\
&&\varepsilon({\hat h}) = \varepsilon({\hat e}) = \varepsilon({\hat f}) =
0,\qquad S({\hat h}) = - {\hat h}, \quad S({\hat e}) = - q^{-1/2}\,{\hat e},
\quad S({\hat f}) = - q^{1/2}\,{\hat f},  \label{eq:qhopf}
\end{eqnarray}
where $[x]_q = (q^x - q^{-x})/(q - q^{-1})$. To fecilitate our later
application, we choose the $(4j+1)$ dimensional irreducible representation
of the ${\cal U}_q(osp(1|2))$ algebra in an asymmetrical manner as
follows: 
\begin{eqnarray}
{\hat h}\;|j\;m> &=& 2 m\;|j\;m>,\qquad {\hat e}\;|j\;m> = |j\;m+1/2>,\qquad 
{\hat e}\;|j\;j> = 0,  \nonumber \\
{\hat f}\;|j\;m> &=&- [j + m]_q\,[[j - m + 1/2]]_q\;|j\;m-1/2>\quad \hbox
{for}\;\;j-m\;\;\hbox {integer},  \nonumber \\
\phantom{{\hat f}\;|j\;m> }&=& [[j + m]]_q\,[j - m + 1/2]_q\;|j\;m-1/2>
\quad \hbox {for}\;\;j-m\;\;\hbox {half-integer},  \label{eq:qrep}
\end{eqnarray}
where $[[x]]_q = (q^{x}-(-1)^{2x}q^{-x})/(q^{1/2}+q^{-1/2})$.

Following the strategy adopted earlier for
constructing the Jordanian deformation of the $sl(N)$ algebra, we give here
the general recipe for obtaining the quantum $R_{\sf h}^{j_1;j_2}$
matrix of an arbitrary representation of the ${\cal U}_{\sf h}(osp(1|2))$
algebra. Explicit demonstration is given for the $1/2 \otimes j$
representation. The relevant $R_{\mathsf{h}}^{1/2;j}$ matrix may be directly
interpreted as the $L$ operator corresponding to the Borel subalgebra of the 
${\cal U}_{\sf h}(osp(1|2))$ algebra. Our construction may,
obviously, be generalized for an arbitrary $j_1 \otimes j_2$ representation.
The primary ingredient for our method is the $R_q^{1/2;j}$ matrix \cite{KR89}
of the ${\cal U}_q(osp(1|2))$ algebra in the $1/2 \otimes j$
representation. A suitable similarity transformation is performed on this $%
R_q^{1/2;j}$ matrix. The transforming matrix is singular in the $q
\rightarrow 1$ limit. \textit{For the transformed matrix, the singularities,
however, systematically cancel yielding a well-defined construction}. The
transformed object, in the $q \rightarrow 1$ limit, directly furnishes the $%
R_{\mathsf{h}}^{1/2;j}$ matrix for the super-Jordanian ${\cal U}_{\sf h}(osp(1|2))$
algebra. Interpreting, as mentioned above, the $R_{\mathsf{h}%
}^{1/2;j}$ matrix obtained here as the $L$ operator corresponding to the
Borel subalgebra of the ${\cal U}_{\sf h}(osp(1|2))$ algebra, we use
the standard FRT procedure \cite{FRT90} to construct the full Hopf structure
of the said Borel subalgebra. The $R_q^{1/2;j}$ matrix of the tensored $1/2
\otimes j$ representation of the ${\cal U}_q(osp(1|2))$ algebra reads 
\cite{KR89} 
\begin{equation}
R_q^{\frac{1}{2};j} = \left(%
\begin{array}{ccc}
q^{\hat h} & -\omega q^{{\hat h}/2}{\hat f} & - \omega \left(1+q^{-1}\right){%
\hat f}^2 \\ 
\noalign{\medskip}0 & 1 & \omega q^{-{(\hat h+1)}/2}{\hat f} \\ 
\noalign{\medskip}0 & 0 & q^{-{\hat h}}%
\end{array}%
\right),  \label{eq:Rqfunarb}
\end{equation}
where $\omega = q - q^{-1}$. We now introduce a transforming matrix $M$,
singular in the $q \rightarrow 1$ limit, as 
\begin{equation}
M = {\sf E}_{q^2}(\eta {\hat e}^2),  \label{eq:singop}
\end{equation}
where $\eta = \frac{\sf h}{q^2 -1}$.  
For any finite value of $j$ the series (49) may be terminated
after setting ${\hat e}^{4j+1} = 0$. As the transforming operator $M$ in (%
\ref{eq:singop}) depends only on the generator ${\hat e}$, our subsequent
results assume \textit{simplest} form for the asymmetric choice of the
representation (\ref{eq:qrep}). Our contraction scheme, however, remains
valid independent of the choice of representation. The $R_q^{j_1;j_2}$
matrix of the ${\cal U}_{\sf h}(osp(1|2))$ algebra may now be subjected to a
similarity transformation followed by a limiting process 
\begin{equation}
\tilde{R}_{\mathsf{h}}^{j_1;j_2} \equiv \lim_{q \to 1} \left[\left(
M_{j_1}^{-1} \otimes M_{j_2}^{-1} \right)\, R_q^{j_1;j_2}\,\left( M_{j_1}
\otimes M_{j_2} \right)\right].  \label{eq:gencon}
\end{equation}
In the followings we will present explicit results for the operator ${\tilde 
R}_{\sf h}^{1/2;j}$. In performing the similarity transformation (\ref{eq:gencon}) we may 
choose any suitable operator ordering. Specifically, starting from left we maintain the 
order ${\hat e}\prec{\hat h}\prec{\hat f}$. In our calculation a class of operators 
${\cal T}_{(\alpha)}=(E_{q^2}(\eta{\hat e}^2))^{-1}\,E_{q^2}(q^{2\alpha}\eta {\hat e}^2)$
satisfying ${\cal T}_{(\alpha+\beta)}q^{(\alpha+\beta){\hat h}}={\cal T}_{(\alpha)}q^{\alpha
{\hat h}}\;{\cal T}_{(\beta)}q^{\beta{\hat h}}$, play an important role. To evaluate $q \to 1$ 
limiting value of the operator ${\cal T}_{(\alpha)}$, we use the identity $E_{q^2}(q^2\,\eta\,
{\hat e}^2)-E_{q^2}(q^{-2}\,\eta\,{\hat e}^2)=\eta\;(q^2-q^{-2})\;{\hat e}^2\;E_{q^2}(\eta
{\hat e}^2)\;\Longrightarrow\;{\cal T}_{(1)}-{\cal T}_{(-1)}=\eta\;(q^2-q^{-2})\;{\hat e}^2$.
Evaluating term by term, the limiting values of ${\cal T}_{(\pm 1)}|_{q\rightarrow 1}\; 
\left(\equiv{\tilde T}_{(\pm 1)}\right)$ are found to be \emph{finite}; and, for these finite 
operators yields ${\tilde T}_{(\pm \alpha)} = \left({\tilde T}_{(\pm 1)}
\right)^{\alpha}$, where ${\tilde T}_{(\alpha)}=\lim_{q \to 1}{\cal T}_{(\alpha)}$.
Writing ${\tilde T}_{(\pm 1)} = {\tilde T}^{\pm 1}$ henceforth, we
immediately observe that in the $q \to 1$ limit, the above identity assumes 
the form 
\begin{equation}
{\tilde T} - {\tilde T}^{-1} = 2{\sf h} e^2\qquad \Longrightarrow \qquad {
\tilde T}^{\pm 1} = \pm{\sf h} e^2 + \sqrt{1 + {\sf h}^2 e^4}.
\label{eq:Tsolve}
\end{equation}
This is our crucial result. Two other operator identities playing key roles
are listed below: 
\begin{eqnarray}
{\hat f}\,{\hat e}^{2n} &=& {\hat e}^{2n}\,{\hat f} - \frac{q}{q + 1}%
\,\{n\}_{q^2}\,{\hat e}^{2n-1}\,{\hat t} - \frac{1}{q + 1}\,\{n\}_{q^{-2}}\,{%
\hat e}^{2n-1}\,{\hat t}^{-1}, \\
{\hat f}^2\,{\hat e}^{2n} &=& {\hat e}^{2n}\,{\hat f}^2 + q \frac{q-1}{q+1}%
\,\{n\}_{q^2}\,{\hat e}^{2n-1}\,{\hat t}\,{\hat f} - q^{-1} \frac{q-1}{q+1}%
\,\{n\}_{q^{-2}}\, {\hat e}^{2n-1}\,{\hat t}^{-1}\,{\hat f}  \nonumber \\
&\phantom{=}& + \frac{q}{q+1}\,\Big(\frac{1}{\omega}\{n\}_{q^4} - q^2\,\frac{%
q-1}{q+1}\,\Big\{
\begin{array}{c}
n \\ 
2%
\end{array}%
\Big\}_{q^2}\Big)\, {\hat e}^{2(n-1)}\,{\hat t}^2  \nonumber \\
&\phantom{=}& - \frac{1}{q+1}\,\Big(\frac{1}{\omega}\{n\}_{q^{-4}} - q^{-2}\,%
\frac{q-1}{q+1}\, \Big\{%
\begin{array}{c}
n \\ 
2%
\end{array}%
\Big\}_{q^{-2}}\Big)\, {\hat e}^{2(n-1)}\,{\hat t}^{-2} 
= - \frac{q}{(q+1)^3}\,\left(q\{n\}_{q^2} +
\{n\}_{q^{-2}}\right)\,{\hat e}^{2(n-1)},  \label{eq:OPE}
\end{eqnarray}
where $\{x\}_q =\frac{1 - q^x}{1 - q}, \{n\}_q! = \{n\}_q\,\{n-1\}_q
\cdots\{1\}_q,\Big\{%
\begin{array}{c}
n \\ 
m%
\end{array}%
\Big\}_q = \frac{\{n\}_q!}{\{n-m\}_q!\{m\}_q!}$ and ${\hat t}^{\pm 1} =
q^{\pm {\hat h}}$. Using the above identities systematically and passing to
the limit $q \to 1$, it may be shown that in our construction of the
operator ${\tilde R}_{\sf h}^{1/2;j}$ \textit{via} (\ref{eq:gencon}), 
\emph{all singularities cancel} yielding in a well-defined answer 
\begin{equation}
{\tilde R}_{\sf h}^{\frac{1}{2};j} =\left(
\begin{array}{ccc}
{\tilde T} & {\sf h} {\tilde T}^{\frac{1}{2}} e & - {\sf h} {\tilde H}
+\frac{\sf h}{4}\left({\tilde T} - {\tilde T}^{-1}\right)\cr 
\noalign{\medskip}0 & 1 & -{\sf h}{\tilde T}^{-\frac{1}{2}} e\cr 
\noalign{\medskip}0 & 0 & {\tilde T}^{-1}
\end{array}
\right),  \label{eq:contR}
\end{equation}
where ${\tilde H} = \frac{1}{2} ({\tilde T} + {\tilde T}^{-1})\;h = \sqrt{1
+ \mathsf{h}^2 e^4}\;h$. One way of interpreting (\ref{eq:contR}) is to
consider it a recipe for obtaining the finite dimensional $R_{\sf h}$
matrices of the ${\cal U}_{\sf h}(osp(1|2))$ algebra. For instance,
using the classical $j=1$ representation, obtained from (\ref{eq:qrep}) in
the $q \rightarrow 1$ limit, we obtain the $R_{\sf h}^{1/2;1} \left(={\tilde 
R}_{\sf h}^{1/2;1}\right)$ matrix as follows: 
\begin{equation}
R_{\sf h}^{\frac{1}{2};1} = \left(%
\begin{array}{ccccccccccccccc}
1 & 0 & \mathsf{h} & 0 & \frac{{\sf h}^2}2 & 0 & \mathsf{h} & 0 & \frac{%
\mathsf{h}^2}2 & 0 & -2\mathsf{h} & 0 & \frac{\mathsf{h}^2}2 & 0 & \mathsf{h}%
^3 \\ 
\noalign{\medskip}0 & 1 & 0 & \mathsf{h} & 0 & 0 & 0 & \mathsf{h} & 0 & 
\frac{\mathsf{h}^2}2 & 0 & -\mathsf{h} & 0 & \frac{\mathsf{h}^2}2 & 0 \\ 
\noalign{\medskip}0 & 0 & 1 & 0 & \mathsf{h} & 0 & 0 & 0 & \mathsf{h} & 0 & 0
& 0 & 0 & 0 & \frac{\mathsf{h}^2}2 \\ 
\noalign{\medskip}0 & 0 & 0 & 1 & 0 & 0 & 0 & 0 & 0 & \mathsf{h} & 0 & 0 & 0
& \mathsf{h} & 0 \\ 
\noalign{\medskip}0 & 0 & 0 & 0 & 1 & 0 & 0 & 0 & 0 & 0 & 0 & 0 & 0 & 0 & 2%
\mathsf{h} \\ 
\noalign{\medskip}0 & 0 & 0 & 0 & 0 & 1 & 0 & 0 & 0 & 0 & 0 & -\mathsf{h} & 0
& \frac{\mathsf{h}^2}2 & 0 \\ 
\noalign{\medskip}0 & 0 & 0 & 0 & 0 & 0 & 1 & 0 & 0 & 0 & 0 & 0 & -\mathsf{h}
& 0 & \frac{\mathsf{h}^2}2 \\ 
\noalign{\medskip}0 & 0 & 0 & 0 & 0 & 0 & 0 & 1 & 0 & 0 & 0 & 0 & 0 & -%
\mathsf{h} & 0 \\ 
\noalign{\medskip}0 & 0 & 0 & 0 & 0 & 0 & 0 & 0 & 1 & 0 & 0 & 0 & 0 & 0 & -%
\mathsf{h} \\ 
\noalign{\medskip}0 & 0 & 0 & 0 & 0 & 0 & 0 & 0 & 0 & 1 & 0 & 0 & 0 & 0 & 0
\\ 
\noalign{\medskip}0 & 0 & 0 & 0 & 0 & 0 & 0 & 0 & 0 & 0 & 1 & 0 & -\mathsf{h}
& 0 & \frac{\mathsf{h}^2}2 \\ 
\noalign{\medskip}0 & 0 & 0 & 0 & 0 & 0 & 0 & 0 & 0 & 0 & 0 & 1 & 0 & -%
\mathsf{h} & 0 \\ 
\noalign{\medskip}0 & 0 & 0 & 0 & 0 & 0 & 0 & 0 & 0 & 0 & 0 & 0 & 1 & 0 & -%
\mathsf{h} \\ 
\noalign{\medskip}0 & 0 & 0 & 0 & 0 & 0 & 0 & 0 & 0 & 0 & 0 & 0 & 0 & 1 & 0
\\ 
\noalign{\medskip}0 & 0 & 0 & 0 & 0 & 0 & 0 & 0 & 0 & 0 & 0 & 0 & 0 & 0 & 1%
\end{array}%
\right).  \label{eq:Rhhalfone}
\end{equation}
The matrix (\ref{eq:contR}) may also be interpreted as the $L$ operator of
the ${\cal U}_{\sf h}(osp(1|2))$ algebra. To this end, we first use
the following invertible map of the quantum ${\cal U}_{\sf h}(osp(1|2))$ algebra 
(\ref{eq:r2hopf}) on the classical algebra (\ref{eq:clalg}): 
\begin{equation}
E = e,\qquad H={\tilde H},\qquad F=f+\frac{\sf h}{4}\,\left(\frac{{\tilde T}-1}
{{\tilde T}+1}\right)\,e - \frac{\sf h}{2}\,\left(\frac{{\tilde T}-1}{{\tilde T}+1}
\right)\,eh,\quad T = {\tilde T},\qquad Y = F^2.  
\label{eq:contrmap}
\end{equation}
The map (\ref{eq:contrmap}) satisfies the algebraic relations (\ref{eq:r2hopf}); and the 
corresponding twist operator may also be determined. Using the map (\ref{eq:contrmap}) the 
operator (\ref{eq:contR}) may be recast in terms of the deformed generators of the 
super-Jordanian ${\cal U}_{\sf h}(osp(1|2))$ algebra as 
\begin{equation}
L \equiv \tilde{R}_{\mathsf{h}}^{\frac{1}{2};j} = \left(%
\begin{array}{ccc}
T & {\sf h} T^{\frac{1}{2}} E & -{\sf h} H +\frac{\sf h}{4}(T-T^{-1})\cr 
\noalign{\medskip}0 & 1 & -{\sf h} T^{-\frac{1}{2}} E \cr 
\noalign{\medskip}0 & 0 & T^{-1}
\end{array}\right).  \label{eq:Lhalf}
\end{equation}
The above $L$ operator of the ${\cal U}_{\sf h}(osp(1|2))$ algebra
has not been obtained before. It allows immediate construction of the full
Hopf structure of the Borel subalgebra of the ${\cal U}_{\sf h}(osp(1|2))$ algebra 
\textit{via} the standard FRT formalism \cite{FRT90}. The algebraic relations for the 
generators ($H, E, T^{\pm 1}$) of the Borel subalgebra is given by 
\begin{equation}
R_{\sf h}^{\frac{1}{2};\frac{1}{2}}\,L_{1}\,L_{2} = L_{2}\,L_{1}\,R_{\sf h}^{\frac{1}{2};
\frac{1}{2}},  \label{eq:RLL}
\end{equation}
where $\ZZ_2$ graded tensor product has been used in defining the operators: $L_{1}=L 
\otimes 1,\, L_{2} =1\otimes L$. The coalgebraic properties of the said Borel subalgebra 
may be succinctly expressed as 
\begin{equation}
\Delta(L) = L \dot{\otimes} L,\qquad \varepsilon(L) = 1,\qquad S(L) =
L^{-1},  \label{eq:Lcoalg}
\end{equation}
where $L^{-1}$ is given by 
\begin{equation}
L^{-1} = \left(%
\begin{array}{ccc}
T^{-1} & -{\sf h} T^{-\frac{1}{2}} E & {\sf h} H + \frac{\sf h}{4}\left(T-T^{-1}\right) \\ 
\noalign{\medskip}0 & 1 & {\sf h} T^{\frac{1}{2}} E \\ 
\noalign{\medskip}0 & 0 & T\end{array}
\right).  \label{eq:Linv}
\end{equation}
This completes our construction of the Hopf structure of the Borel
subalgebra of the super-Jordanian ${\cal U}_{\sf h}(osp(1|2))$
algebra, obtained by deforming the $r_2$ matrix, by employing the
contraction scheme described earlier. Our results fully coincide with the
Hopf structure given in (\ref{eq:r2hopf}). This validates our contraction
scheme elaborated earlier. Our recipe (\ref{eq:gencon}) for obtaining the 
$R_{\sf h}^{j_1;j_2}$ matrix for a $j_1 \otimes j_2$ representation of
the super-Jordanian ${\cal U}_{\sf h}(osp(1|2))$ algebra may be
continued arbitrarily. The matrices such as $R_{\sf h}^{1;j}$ may be
interpreted as higher dimensional $L$ operators \cite{JJ96} obeying duality
relations with relevant ${\cal T}$ matrices.

\subsection{Jordanian ${\bf U_h}(osp(1|2))$ algebra: a nonlinear realization}

The classical $r_1$ matrix has been quantized earlier \cite{CK98} using the
inclusion $sl(2) \subset osp(1|2)$. These authors have expressed the
resultant triangular deformed $osp(1|2)$ superalgebra in terms of the
classical basis set. On the other hand Ohn \cite{P4} employed a nonlinear
basis set to formulate the deformed Jordanian ${\cal U}_{\sf h}(sl(2))$ algebra. Consequently, 
Ohn's ${\cal U}_{\sf h}(sl(2))$ algebra \cite{P4} do not directly appear as a Hopf subalgebra 
of the deformed $osp(1|2)$ superalgebra considered in Ref. \cite{CK98}. Moreover, if the 
algebraic relations are described in terms of the undeformed classical basis set, the coproduct 
structures tend to be complicated in nature.

In the following, we present the quantized Hopf structure corresponding to
the classical $r_1$ matrix in terms of {\it nonlinear basis elements}.
Our algebra {\it explicitly includes} Ohn's ${\cal U}_{\sf h}(sl(2))$ algebra as a Hopf 
subalgebra. The coproduct structure we obtain is \textit{considerably simple}. To distinguish 
the Jordanian deformed $osp(1|2)$ superalgebra, its generators and the deformation parameter 
considered here from the corresponding objects displayed above, we express them in boldfaced 
notations. The Jordanian ${\bf U_h}(osp(1|2))$ algebra corresponding to the $r_1$ matrix is 
generated by the elements $({\bf H, E, F, X, Y})$. Their classical analogs are $(h, e, f, 
b_{\pm})$ respectively. The elements ${\bf T}^{\pm 1} = \exp ({\pm {\bf hX}})$
are also introduced. The deformation parameter is denoted by ${\bf h}$.
The Hopf structure of the ${\bf U_h}(osp(1|2))$ algebra is obtained by
maintaining the following properties: ({\sf i}) In the classical limit
the quantum coproduct map conforms to the classical cocommutator. ({\sf ii}) 
The coproduct map is a homomorphism of the algebra, and it satisfies the
coassociativity constraint. ({\sf iii}) Generator ${\bf X}$ is the only primitive element. 
The commutation relations of the ${\bf U_h}(osp(1|2))$ algebra reads [23] 
\begin{eqnarray}
&& [\mathbf{H},\,\mathbf{E}] = \frac {1}{2} (\mathbf{T}+\mathbf{T}^{-1})%
\mathbf{E},  \nonumber \\
&& [\mathbf{H},\,\mathbf{F}] = -\frac {1}{4} (\mathbf{T}+\mathbf{T}^{-1})%
\mathbf{F} -\frac {1}{4}\mathbf{F} (\mathbf{T}+\mathbf{T}^{-1}) -\frac{%
\mathbf{h}}{8} ((\mathbf{T}-\mathbf{T}^{-1}) \mathbf{H} + \mathbf{H} (%
\mathbf{T}-\mathbf{T}^{-1})) \mathbf{E}  \nonumber \\
&&\phantom{[{\bf H},\,{\bf F}]=} - \frac{\mathbf{h}}{8} \mathbf{E} ((\mathbf{%
T}-\mathbf{T}^{-1}) \mathbf{H} + \mathbf{H}(\mathbf{T}-\mathbf{T}^{-1})), 
\nonumber \\
&& \{\mathbf{E},\,\mathbf{F}\} = -\frac {1}{4} (\mathbf{T}+\mathbf{T}^{-1}) 
\mathbf{H} - \frac {1}{4} \mathbf{H} (\mathbf{T}+\mathbf{T}^{-1}),\qquad [%
\mathbf{H},\,\mathbf{T}^{\pm 1}] = \mathbf{T}^{\pm 2} - 1,  \nonumber \\
&&[\mathbf{H},\,\mathbf{Y}] = - \frac{1}{2} (\mathbf{T} + \mathbf{T}^{-1}) 
\mathbf{Y} - \frac{1}{2} \mathbf{Y} (\mathbf{T} + \mathbf{T}^{-1}),\qquad [%
\mathbf{T}^{\pm 1},\,\mathbf{Y}] = \pm \frac{\mathbf{h}}{2} (\mathbf{T}^{\pm
1} \mathbf{H} + \mathbf{H } \mathbf{T}^{\pm 1}),  \nonumber \\
&&\mathbf{E}^2=\frac {1}{2\mathbf{h}} (\mathbf{T}-\mathbf{T}^{-1}), \qquad [%
\mathbf{Y}, \mathbf{E}] = \mathbf{F},\qquad [\mathbf{T}^{\pm 1}, \mathbf{F}]
= \pm \frac{\mathbf{h}}{2} (\mathbf{T}^{\pm 2} + 1) \mathbf{E},  \nonumber \\
&& \mathbf{F}^2 = -\mathbf{Y}+\frac{\mathbf{h}}{8} (\mathbf{T}-\mathbf{T}%
^{-1}) \mathbf{H}^2 + \frac{\mathbf{h}}{4} (\mathbf{T}-\mathbf{T}^{-1}) 
\mathbf{E F} + \frac{3\mathbf{h}}{16}(\mathbf{T}^2-\mathbf{T}^{-2}) \mathbf{H%
} + \frac{\mathbf{h}}{4} (\mathbf{T}-\mathbf{T}^{-1})  \nonumber \\
&&\phantom{{\bf F}^2 =} + \frac{9\mathbf{h}}{128} (\mathbf{T}-\mathbf{T}%
^{-1})^3,  \nonumber \\
&&[\mathbf{F},\;\mathbf{Y}] = \frac{\mathbf{h}}{4} (\mathbf{T}-\mathbf{T}%
^{-1}) \mathbf{F} + \frac{\mathbf{h}}{2} (\mathbf{T}-\mathbf{T}^{-1}) 
\mathbf{E Y} - \frac{\mathbf{h}^2}{4}\mathbf{E H}^2 - \frac{3\mathbf{h}^2}{8}
(\mathbf{T}+\mathbf{T}^{-1}) \mathbf{E H} - \frac{\mathbf{h}^2}{2} \mathbf{E}
\nonumber \\
&&\phantom{[{\bf F},\;{\bf Y}] =} - \frac{15\mathbf{h}^2}{64} (\mathbf{T}-%
\mathbf{T}^{-1})^2 \mathbf{E}  \label{eq:r1alg}
\end{eqnarray}
and the corresponding coalgebraic structure is given by [23]
\begin{eqnarray}
&&\Delta(\mathbf{H}) = \mathbf{H} \otimes \mathbf{T} + \mathbf{T}^{-1}
\otimes \mathbf{H},\qquad\qquad\qquad\qquad \Delta(\mathbf{E}) = \mathbf{E}
\otimes \mathbf{T}^{-1/2} + \mathbf{T}^{1/2}\otimes \mathbf{E},  \nonumber \\
&&\Delta (\mathbf{F}) = \mathbf{F} \otimes \mathbf{T}^{1/2} + \mathbf{T}%
^{-1/2}\otimes \mathbf{F}+ \frac{\mathbf{h}}{4} \mathbf{T}^{-1} \mathbf{E}
\otimes \left(\mathbf{T}^{-1/2}\mathbf{H} + \mathbf{H T}^{-1/2}\right) - 
\frac{\mathbf{h}}{4} \left (\mathbf{T}^{1/2} \mathbf{H} + \mathbf{H T}%
^{1/2}\right)\otimes \mathbf{T E},  \nonumber \\
&&\Delta(\mathbf{T}^{\pm 1}) = \mathbf{T}^{\pm 1} \otimes \mathbf{T}^{\pm
1}, \qquad\qquad\qquad\qquad \Delta(\mathbf{Y}) = \mathbf{Y} \otimes \mathbf{%
T} + \mathbf{T}^{-1} \otimes \mathbf{Y},  \nonumber \\
&&\varepsilon(\mathbf{H}) = \varepsilon(\mathbf{E}) = \varepsilon(\mathbf{F}%
) = \varepsilon(\mathbf{Y}) = 0,\qquad\qquad\qquad \varepsilon(\mathbf{T}%
^{\pm 1}) = 1,  \nonumber \\
&&S(\mathbf{H}) = - \mathbf{H} + 2 \mathbf{h} \mathbf{E}^2,\qquad\qquad S(%
\mathbf{E}) = - \mathbf{E},\qquad\qquad S(\mathbf{F}) = - \mathbf{F} - \frac{%
\mathbf{h}}{2}\left(\mathbf{T} + \mathbf{T}^{-1}\right) \mathbf{E}, 
\nonumber \\
&&S(\mathbf{T}^{\pm 1}) = \mathbf{T}^{\mp 1},\qquad\qquad\qquad\qquad S(%
\mathbf{Y}) = - \mathbf{Y} - \mathbf{h H} + \mathbf{h}^2 \mathbf{E}^2.
\label{eq:r1coalg}
\end{eqnarray}
All the Hopf algebra axioms can be verified by direct calculation. The
universal ${\bf R_h}$ matrix of the Jordanian ${\bf U_h}(osp(1|2))$ 
algebra is of the factorized form [26]: 
\begin{equation}
\mathbf{R_h} = G_{21}^{-1} G,\qquad G = \exp\left(\mathbf{h}\,\mathbf{T H}
\otimes \mathbf{X}\right),  \label{eq:r1univR}
\end{equation}
which coincides with the universal ${\cal R}_{\sf h}$ matrix of the ${\cal U}_{\sf h}
(sl(2))$ subalgebra [26] involving the highest weight root vector.

Before discussing the general structure of a \textit{class} of invertible
maps of the ${\bf U_h}(osp(1|2))$ algebra on the classical ${\cal U}
(osp(1|2))$ algebra, we notice that the comultiplication map of a set of
three operators 
\begin{equation}
\mathbf{T}^{-1/2} \mathbf{E},\quad \mathbf{T H}, \quad \mathbf{T}^{1/2} 
\mathbf{F} + \frac{\mathbf{h}}{8} \mathbf{T}^{1/2} (\mathbf{T} - \mathbf{T}%
^{-1}) \mathbf{E} - \frac{\mathbf{h}}{2} \mathbf{T}^{1/2} \mathbf{E H},
\label{eq:mtmop}
\end{equation}
when acted by the twist operator corresponding to the factorized form of the
universal $\mathbf{R_h}$ matrix given in (\ref{eq:r1univR}), reduce to the
classical cocommutative coproduct: 
\begin{equation}
G\Delta({\cal X})G^{-1}={\cal X}\otimes\mathbf{I}+\mathbf{I}\otimes{\cal X},\label{eq:mtmtwist}
\end{equation}
where ${\cal X}$ is an element of the set (\ref{eq:mtmop}). From the
commutation rules (\ref{eq:r1alg}) it becomes evident that the operators in
the set (\ref{eq:mtmop}) satisfy the classical algebra generated by $(e, h,
f)$ respectively. These operators constitute an important special case of a
class of maps discussed below.

To construct a \textit{class} of maps interrelating the Jordanian ${\bf U_h}(osp(1|2))$ 
algebra, obtained \textit{via} quantization of the $r_1$-matrix, and its classical analog 
${\cal U}(osp(1|2))$ we proceed {\it via} an ansatz as follows: 
\begin{eqnarray}
&&\mathbf{E} = \varphi_1 (b_{+}) e,\qquad \mathbf{H} = \varphi_2(b_{+}) h, 
\qquad
\mathbf{F} = \varphi_3(b_{+}) f + u_1(b_{+}) e + u_2(b_{+}) e h,
\label{eq:ansatz}
\end{eqnarray}
where the `mapping functions' ($\varphi_1, \varphi_2, \varphi_3; u_1, u_2$)
depend only on the classical generator $b_{+}$. In the classical limit 
$\mathbf{h} \to 0$ the above functions satisfy the property: ($\varphi_1,
\varphi_2, \varphi_3; u_1, u_2$) $\to$ ($1, 1, 1; 0, 0$). The operators 
$\mathbf{T}^{\pm 1}$ may now be expressed as 
\begin{equation}
\mathbf{T}^{\pm 1} = \pm \mathbf{h} b_{+} (\varphi_1(b_{+}))^2 + \sqrt{ 1 + 
\mathbf{h}^2 b_{+}^2 (\varphi_1(b_{+}))^4}.  \label{eq:TiTin}
\end{equation}
Substituting the ansatz (\ref{eq:ansatz}) in the defining relations (\ref%
{eq:r1alg}) for the $\mathbf{U_h}(osp(1|2))$ algebra we, for a \textit{given}
function $\varphi_1$, obtain a set of \textit{six} nonlinear equations for 
\textit{four} unknown functions: 
\begin{eqnarray}
&&\left(\varphi_1(b_{+}) + 2 b_{+} \varphi_1^{\prime}(b_{+})\right)
\varphi_2(b_{+}) - \sqrt{ 1 + \mathbf{h}^2 b_{+}^2 (\varphi_1(b_{+}))^4}\,
\varphi_1(b_{+}) = 0,  \nonumber \\
&&2 b_{+} \varphi_2(b_{+}) \varphi_3^{\prime}(b_{+}) - \varphi_2(b_{+})
\varphi_3(b_{+}) + \sqrt{ 1 + \mathbf{h}^2 b_{+}^2 (\varphi_1(b_{+}))^4}\,
\varphi_3(b_{+})= 0,  \nonumber \\
&&2 b_{+} \varphi_2(b_{+}) u_1^{\prime}(b_{+}) + \left(\varphi_2(b_{+}) + 
\sqrt{1 + \mathbf{h}^2 b_{+}^2 (\varphi_2(b_{+}))^4}\right)\,u_1(b_{+}) 
\nonumber \\
&&\qquad\qquad + \mathbf{h}^2 b_{+} \sqrt {1 + \mathbf{h}^2 b_{+}^2
(\varphi_1(b_{+}))^4}\, (\varphi_1(b_{+}))^3 = 0,  \nonumber \\
&&\left(\varphi_2(b_{+}) - 2 b_{+}\varphi_2^{\prime}(b_{+}) + \sqrt {1 + 
\mathbf{h}^2 b_{+}^2 (\varphi_1(b_{+}))^4} \right) u_2(b_{+}) +
\varphi_2^{\prime}(b_{+}) \varphi_3(b_{+})  \nonumber \\
&&\qquad\qquad + 2 b_{+} \varphi_2(b_{+}) u_2^{\prime}(b_{+}) + \mathbf{h}^2
b_{+} (\varphi_1(b_{+}))^3 \varphi_2(b_{+}) = 0,  \nonumber \\
&&\varphi_1(b_{+}) \left(2 b_{+} u_2(b_{+}) - \varphi_3(b_{+})\right) + 
\sqrt {1 + \mathbf{h}^2 b_{+}^2 (\varphi_1(b_{+}))^4}\, \varphi_2(b_{+}) = 0,
\nonumber \\
&&b_{+} \varphi_1(b_{+}) \left(2 u_1(b_{+}) + u_2(b_{+})\right) - b_{+}
\varphi_1^{\prime}(b_{+}) \left( \varphi_3(b_{+}) - 2 b_{+}
u_2(b_{+})\right) + \mathbf{h}^2 b_{+}^2 (\varphi_1(b_{+}))^4 = 0.
\label{eq:fordiff}
\end{eqnarray}
Maintaining the classical limit the above set of equations may be
consistently solved as follows: 
\begin{eqnarray}
&&\varphi_2(b_{+}) = \frac{\sqrt{ 1 + \mathbf{h}^2 b_{+}^2
(\varphi_1(b_{+}))^4}\, \varphi_1(b_{+})} {\varphi_1(b_{+}) + 2 b_{+}
\varphi_1^{\prime}(b_{+})},\qquad \varphi_3(b_{+}) = \frac{1}{%
\varphi_1(b_{+})},  \nonumber \\
&&u_1(b_{+}) = - \frac{\mathbf{h}^2}{4}\,b_{+} (\varphi_1(b_{+}))^3,\qquad
u_2(b_{+}) = \frac{1 - \sqrt{ 1 + \mathbf{h}^2 b_{+}^2 (\varphi_1(b_{+}))^4}%
\, \varphi_2(b_{+})}{2 b_{+} \varphi_1(b_{+})}.  \label{eq:forsol}
\end{eqnarray}
Inverse maps expressing the classical generators in terms of the relevent
quantum generators are obtained by assuming the ansatz 
\begin{eqnarray}
&&e = \psi_1(\mathbf{T})\,\mathbf{E},\qquad h = \psi_2(\mathbf{T})\,\mathbf{H},  \qquad
f=\psi_3(\mathbf{T})\,\mathbf{F}+w_1(\mathbf{T})\,\mathbf{E}+w_2(\mathbf{T})\,\mathbf{E H},
\label{eq:invan}
\end{eqnarray}
where $(\psi_1, \psi_2, \psi_3; w_1, w_2)$ are functions of $\mathbf{T}$
obeying the limiting property: $(\psi_1, \psi_2, \psi_3; w_1, w_2) \to (1,
1, 1; 0, 0)$ as $\mathbf{h} \to 0$. The differential equations obeyed by the
`mapping functions' introduced in (\ref{eq:invan}) read 
\begin{eqnarray}
&&2 (\mathbf{T}^2 - 1)\,\psi_1^{\prime}(\mathbf{T}) \psi_2(\mathbf{T}) + (%
\mathbf{T} + \mathbf{T}^{-1})\,\psi_1(\mathbf{T}) \psi_2(\mathbf{T}) - 2
\psi_1(\mathbf{T}) = 0,  \nonumber \\
&&2 (\mathbf{T}^2 - 1)\,\psi_2(\mathbf{T}) \psi_3^{\prime}(\mathbf{T}) - (%
\mathbf{T} + \mathbf{T}^{-1})\,\psi_2(\mathbf{T}) \psi_3(\mathbf{T}) +2
\psi_3(\mathbf{T}) = 0,  \nonumber \\
&&\psi_2(\mathbf{T})\,\left(2 (\mathbf{T} + \mathbf{T}^{-1})\,w_1(\mathbf{T}%
) + 4 (\mathbf{T}^2 - 1)\,w_1^{\prime}(\mathbf{T}) - \mathbf{h} (\mathbf{T}%
^2 - \mathbf{T}^{-2})\,\psi_3(\mathbf{T})\right) + 4 w_1(\mathbf{T}) = 0, 
\nonumber \\
&&(\mathbf{T} - \mathbf{T}^{-1})\,\psi_2(\mathbf{T})\, \left(2\mathbf{T}%
\,w_2^{\prime}(\mathbf{T}) - \mathbf{h}\,\psi_3(\mathbf{T})\right) + \mathbf{%
h} (\mathbf{T}^2 + 1)\,\psi_2^{\prime}(\mathbf{T})\, \psi_3(\mathbf{T}) 
\nonumber \\
&&\qquad\qquad + \left((\mathbf{T} + \mathbf{T}^{-1})\,\psi_2(\mathbf{T}) -
2 (\mathbf{T}^2 - 1) \,\psi_2^{\prime}(\mathbf{T}) + 2\right) w_2(\mathbf{T}%
) = 0,  \nonumber \\
&&2 (\mathbf{T} - \mathbf{T}^{-1})\,\psi_1(\mathbf{T}) w_2(\mathbf{T}) - 
\mathbf{h} (\mathbf{T} + \mathbf{T}^{-1})\,\psi_1(\mathbf{T})\,\psi_3(%
\mathbf{T}) + 2 \mathbf{h} \psi_2(\mathbf{T}) = 0,  \nonumber \\
&&(\mathbf{T} - \mathbf{T}^{-1})\,\psi_1(\mathbf{T})\,\left(4 w_1(\mathbf{T}%
) + (\mathbf{T} + \mathbf{T}^{-1})\,w_2(\mathbf{T}) - \mathbf{h} (\mathbf{T}
- \mathbf{T}^{-1})\,\psi_3(\mathbf{T})\right)  \nonumber \\
&&\qquad\qquad + \mathbf{T} (\mathbf{T} - \mathbf{T}^{-1})\,\psi_1^{\prime}(%
\mathbf{T}) \left(2 (\mathbf{T} - \mathbf{T}^{-1})\,w_2(\mathbf{T}) - 
\mathbf{h} (\mathbf{T} + \mathbf{T}^{-1})\,\psi_3(\mathbf{T})\right) = 0.
\label{eq:backdiff}
\end{eqnarray}
Treating the function $\psi_1(\mathbf{T})$ as known, and maintaining the
limiting properties the remaining functions may be solved uniquely: 
\begin{eqnarray}
\psi_2(\mathbf{T}) &=& \frac{2 \psi_1(\mathbf{T})} {(\mathbf{T} + \mathbf{T}%
^{-1})\,\psi_1(\mathbf{T}) + 2 (\mathbf{T}^2 - 1)\,\psi_1^{\prime}(\mathbf{T}%
)},\qquad\qquad \psi_3(\mathbf{T}) = \frac{1}{\psi_1(\mathbf{T})},  \nonumber
\\
w_1(\mathbf{T}) &=& \frac{\mathbf{h} (\mathbf{T} - \mathbf{T}^{-1})} {8
\psi_1(\mathbf{T})},\qquad\qquad\quad w_2(\mathbf{T}) = \mathbf{h}\,\frac{%
\mathbf{T} + \mathbf{T}^{-1} - 2 \psi_2(\mathbf{T})}{2 (\mathbf{T} - \mathbf{%
T}^{-1})\,\psi_1(\mathbf{T})}.  \label{eq:backsol}
\end{eqnarray}

The general structure of the twisting elements corresponding to the given maps may be 
described as follows. Let $m$ be a deformation map and $m^{-1}$ be its inverse: 
$m: (\mathbf{E, H, F}) \to (e, h, f)$ and $m^{-1}: (e, h, f) \to (\mathbf{E,
H, F})$. The classical $(\Delta_0)$ cocommutative and the quantum $(\Delta)$
non-cocommutative coproducts are related by the twisting element as 
\begin{equation}
{\cal G}\,\Delta \circ m^{-1} (\phi)\,{\cal G}^{-1} = (m^{-1} \otimes
m^{-1}) \circ \Delta_0(\phi)\qquad \forall \phi \in {\cal U}(osp(1|2)),
\label{eq:Drinfeld}
\end{equation}
where the twisting element ${\cal G} \in {\bf U}_h(osp(1|2))^{\otimes
2}$ satisfies the cocycle condition 
\begin{equation}
({\cal G} \otimes{\bf I})\,((\Delta\otimes\hbox {id}){\cal G})=({\bf I}\otimes
{\cal G})\,((\hbox{id} \otimes \Delta){\cal G}).  \label{eq:cocycle}
\end{equation}
Similarly the classical $(S_0)$ and the quantum $(S)$ antipode maps are
related as follows: 
\begin{equation}
{\bf g}\,S \circ m^{-1} (\phi)\,{\bf g}^{-1} = m^{-1} \circ S_0(\phi),\qquad {\bf g} 
\in {\bf U}_h(osp(1|2)).  \label{eq:Srel}
\end{equation}
The transforming operator ${\sf g}$ for the antipode map may be expressed
in terms of the twist operator ${\cal G}$ as 
\begin{equation}
{\bf g} = \mu \circ(\hbox{id} \otimes S) {\cal G},  \label{eq:antwi}
\end{equation}
where $\mu$ is the multiplication map.

The first example belonging to the class of invertible maps discussed
earlier plays a key role in the present construction. With the choice 
\begin{equation}
\varphi_1(b_{+}) = (1 - 2 {\bf h} b_{+})^{-1/4},\qquad \psi_1({\bf T})
= {\bf T}^{-1/2},  \label{eq:mtmchoice}
\end{equation}
we obtain, after implementing (\ref{eq:ansatz}, \ref{eq:forsol}, \ref{eq:invan}) and 
(\ref{eq:backsol}), the following direct map 
\begin{eqnarray}
{\bf E} &=& (1 - 2 {\bf h} b_{+})^{-1/4}\,e,\quad {\bf H}=\sqrt{(1-2{\bf h}b_{+})}\,h,
\quad {\bf T}^{\pm 1}=(1-2{\bf h}b_{+})^{\mp 1/2},  \nonumber \\
{\bf F} &=& (1 - 2 {\bf h} b_{+})^{1/4}\,f-\frac{{\bf h}^2}{4}b_{+} (1 -2{\bf h} 
b_{+})^{-3/4}\,e + \frac{\bf h}{2}\,(1-2{\bf h} b_{+})^{1/4}\,eh  \label{eq:mtmdir}
\end{eqnarray}
and its inverse 
\begin{equation}
e = {\bf T}^{-1/2} {\bf E},\quad h={\bf T H}, \quad f={\bf T}^{1/2}{\bf F} + 
\frac{\bf h}{8}{\bf T}^{1/2}({\bf T}-{\bf T}^{-1}){\bf E}-\frac{\bf h}{2}{\bf T}^{1/2}
{\bf E H}.  \label{eq:mtminv}
\end{equation}
It may be observed from (\ref{eq:mtmtwist}) that the operator $G$
corresponding to the factorized form of the universal ${\bf R_h}$ matrix
given in (\ref{eq:r1univR}) plays the role of the twist operator ${\cal G}$ for the map 
(\ref{eq:mtmdir}) and its inverse. In this sense we refer to it as the `minimal twist map'. 
The operator ${\bf g}$ transforming the antipode map may be, \textit{{\`a} la} 
(\ref{eq:Srel}), explicitly evaluated in a \textit{closed} form: 
\begin{equation}
{\bf g}=\exp\left(-\frac{1}{2}{\bf TH} (1-{\bf T}^{-2})\right).
\label{eq:gvalue}
\end{equation}
Combining (\ref{eq:gvalue}) with the property (\ref{eq:antwi}) we now
immediately obtain a disentanglement relation, which, if expressed in terms
of classical generators, reads as follows: 
\begin{equation}
\mu \left[\exp\left(\frac{1}{2}\,h \otimes \ln (1 - 2{\bf h} b_{+})\right)
\right] = \exp(-{\bf h} h b_{+}).  \label{eq:disentan}
\end{equation}
In the above relation ${\bf h}$ may be treated as an arbitrary parameter.
To our knowledge the above disentanglement formula involving the classical 
$sl(2)$ generators was not observed before.

Another map, where the Cartan element ${\bf H}$ of the Jordanian ${\bf U_h}(osp(1|2))$ 
algebra remains diagonal, is given by `mapping functions' 
\begin{equation}
\varphi_1(b_{+}) = \frac{1}{\sqrt{1 - \frac{\mathbf{h}^2 b_{+}^2}{4}}},
\qquad \psi_1(\mathbf{T}) = \hbox{sech}(\frac{\mathbf{hX}}{2}).
\label{eq:2ndcho}
\end{equation}
Proceeding as in the previous example, we now map the $\mathbf{U_h}(osp(1|2))
$ algebra on the the classical ${\cal U}(osp(1|2))$ algebra as 
\begin{eqnarray}
\mathbf{E} &=& \frac{1}{\sqrt{1 - \frac{\mathbf{h}^2 b_{+}^2}{4}}}\,e,\qquad 
\mathbf{H} = h,\qquad \mathbf{T}^{\pm 1} = \frac{1 \pm \frac {\mathbf{h}
b_{+}}{2}}{1 \mp \frac {\mathbf{h} b_{+}}{2}},  \nonumber \\
\mathbf{F} &=& \sqrt{1 - \frac{\mathbf{h}^2 b_{+}^2}{4}}\,f - \frac{\mathbf{h%
}^2 b_{+}} {4 \left(1 - \frac{\mathbf{h}^2 b_{+}^2}{4}\right)^{\frac{3}{2}}}%
\,e - \frac{\mathbf{h}^2 b_{+}} {4\,\sqrt{1 - \frac{\mathbf{h}^2 b_{+}^2}{4}}%
}\,eh.  \label{eq:hdiadir}
\end{eqnarray}
The inverse map now reads 
\begin{eqnarray}
&&e = \hbox{sech}(\frac{\mathbf{hX}}{2}) \mathbf{E},\qquad h = \mathbf{%
H},\qquad
f = \hbox{cosh}(\frac{\mathbf{hX}}{2})\,\mathbf{F} + \frac{\mathbf{h}}{4} %
\hbox{sinh}(\mathbf{hX})\, \hbox{cosh}(\frac{\mathbf{hX}}{2})\,\mathbf{E} + 
\frac{\mathbf{h}}{2} \hbox{sinh}(\frac{\mathbf{hX}}{2})\,\mathbf{EH}.
\label{eq:hdiainv}
\end{eqnarray}
The twist operator ${\cal G}$ for the map (\ref{eq:hdiadir}), unlike the
previous example of closed-form expression for the `minimal twist map' given
in (\ref{eq:mtmtwist}), may be determined only in a series: 
\begin{equation}
\mathcal{G} = \mathbf{I} \otimes \mathbf{I} + \frac{\mathbf{h}}{2} \mathbf{r}
+ \frac{\mathbf{h}^2}{8} (\mathbf{r}^2 + \mathbf{H} \otimes \mathbf{X}^2 +%
\mathbf{X}^2 \otimes \mathbf{H}) + O(\mathbf{h}^3),  \label{eq:hdiatwi}
\end{equation}
where $\mathbf{r} = \mathbf{H} \otimes \mathbf{X} - \mathbf{X} \otimes 
\mathbf{H}$. The corresponding transforming operator for the antipode map is
determined from (\ref{eq:Srel}): 
\begin{equation}
\mathbf{g} = {} - \mathbf{h X} + \frac{1}{2} \mathbf{h}^2 \mathbf{X}^2 + O(%
\mathbf{h}^3).  \label{eq:hdiant}
\end{equation}
Other suitable maps belonging to the \textit{class} described before may be
obtained by just choosing sufficiently well-behaved, but otherwise \textit{%
arbitrary}, `mapping functions' $\varphi_1(b_{+})$ and its inverse $\psi_1(T)
$. Our formalism given in (\ref{eq:forsol}) and (\ref{eq:backsol}) readily
produces solutions for the other `mapping functions' listed in (\ref%
{eq:ansatz}) and (\ref{eq:invan}). These maps may be used to immediately
generate representations of the Jordanian ${\bf U_h}(osp(1|2))$ algebra.

\section{The Nonstandard superalgebra Enveloping ${\cal U}_{\sf h}(sl(2|1))$}   

The main object of this section is to extend our construction to ${\cal U}(sl(2|1))$ 
envelopping superalgebras. Let us just recall the more important points concerning $sl(2|1)$: 
The superalgebra $sl(1|2)$ is generated by six generators $\{h_1,\;h_2,\;h_3,\;e_1,\;f_1,
\;e_2,\;f_2,\;e_3,\;f_3\}$ and the commutation relations    
\begin{equation}
\begin{array}{lllll}
\left[h_1,\;h_2\right]=0, & & & & (a)\\
\left[h_1,\;h_3\right]=0, &\left[h_2,\;h_3\right]=0, & & &\\
\left[h_1,\;e_1\right]=2e_1,\qquad &\left[h_1,\;f_1\right]=-2f_1,\qquad &
\left[h_1,\;e_2\right]=-e_2,\qquad  &\left[h_1,\;f_2\right]=f_2,\qquad & (b)\\
\left[h_1,\;e_3\right]=e_3, & \left[h_1,\;f_3\right]=-f_3,&&&\\ 
\left[h_2,\;e_1\right]=-e_1,& \left[h_2,\;f_1\right]=f_1,&\left[h_2,\;e_2\right]=0, & 
\left[h_2,\;f_2\right]=0,& (c)\\
\left[h_2,\;e_3\right]=-e_3,& \left[h_2,\;f_3\right]=f_3,&&&\\ 
\left[h_3,\;e_1\right]=e_1,& \left[h_3,\;f_1\right]=-f_1,&\left[h_3,\;e_2\right]=-e_2, &
\left[h_3,\;f_2\right]=f_2&\\
\left[h_3,\;e_3\right]=0,& \left[h_3,\;f_3\right]=0,&&&\\ 
\left[e_1,\;f_1\right]=h_1,& \left[e_2,\;f_2\right]=h_2,& &&(d)\\
\left[e_3,\;f_3\right]=h_3,& &\\
\left[e_1,\;f_2\right]=\left[e_2,\;f_1\right]=0,\qquad & e_2^2=f_2^2=0,&&& (e)\\
\left[e_1,\;e_2\right]=e_3,& \left[f_2,\;f_1\right]=f_3,&&& \\
e_3^2=f_3^2=0,& \left[e_1,\;e_3\right]=0, &\left[f_3,\;f_1\right]=0,& \left[e_2,\;e_3\right]=
0, &\left[f_2,\;f_3\right]=0,\\ 
\left[f_1,\;e_3\right]=e_2,& \left[f_3,\;e_1\right]=f_2, &\left[f_2,\;e_3\right]=e_1, &
\left[f_3,\;e_2\right]=f_1, & 
\end{array}
\end{equation} 
where the commutator $[\;,\;]$ is understood as the $\ZZ_2$-graded one: $[a,\;b]=ab-\left(-
\right)^{\deg(a)\deg(b)}ba$.
The generators $h_1$, $h_2$, $h_3$, $e_1$ and $f_1$ are even ($\deg(h_1)=\deg(h_2)=\deg(h_3)=
\deg(e_1)=\deg(f_1)=0$), while $e_2$, $f_2$, $e_3$ and $f_3$ are odd, ($\deg(e_2)=\deg(f_2)=
\deg(e_3)=\deg(f_3)=1$). As a Hopf superalgebra, the universal enveloping ${\cal U}(sl(2|1))$ 
of $sl(2|1)$ is generated just by six elements: it is sufficient to start from $\{h_1,\;h_2,\;
e_1,\;f_1,\;e_2,\;f_3\}$ restricted by relations (a), (b), (c), (d), (e), (f), (g) only, and 
\begin{eqnarray}
&&e_1^2e_2-2e_1e_2e_1+e_2e_1^2=0, \qquad\qquad 
f_1^2f_2-2f_1f_2f_1+f_2f_1^2=0.
\end{eqnarray}
The two last equations are called the Serre relations. Let us just mention that there is a 
$\CC$-algebra automorphism $\phi$ of ${\cal U}(sl(2|1))$ such that
\begin{equation}
\begin{array}{lll}
\phi\left(e_1\right)=e_1,\qquad\qquad &\phi\left(f_1\right)=f_1,\qquad \qquad 
&\phi\left(h_1\right)=h_1, \\
\phi\left(e_2\right)=f_3,\qquad\qquad &\phi\left(f_2\right)=-e_3,\qquad \qquad 
&\phi\left(h_2\right)=-h_3, \\
\phi\left(e_3\right)=-f_2,\qquad\qquad &\phi\left(f_3\right)=e_2,\qquad \qquad 
&\phi\left(h_3\right)=-h_2.
\end{array}
\end{equation}
The quasitriangular quantum Hopf superalgebra ${\cal U}_q(sl(2|1))$ ($q$ is an arbitrary 
complex number), by analogy with ${\cal U}(sl(2|1))$, is generated by six elements 
$\{{\hat h}_1,\;{\hat h}_2,\;{\hat e}_1,\;{\hat e}_2,\;{\hat f}_1,\;{\hat f}_2\}$ under the 
relations 
\begin{equation}
\begin{array}{ll}
[{\hat h}_1,\;{\hat h}_2]=0,&\phantom{xx}\\
\noalign{\medskip} [{\hat h}_1,\;{\hat e}_1]=2{\hat e}_1,\qquad\qquad\qquad & 
[{\hat h}_1,\;{\hat f}_1]=-2{\hat f}_1,\\
\noalign{\medskip} [{\hat h}_1,\; {\hat e}_2]=-{\hat e}_2, & 
[{\hat h}_1,\;{\hat f}_2]={\hat f}_2,\\
\noalign{\medskip} [{\hat h}_2,\; {\hat e}_1]=-{\hat e}_1, & 
[{\hat h}_2,\;{\hat f}_1]={\hat f}_1,\\
\noalign{\medskip} [{\hat h}_2,\; {\hat e}_2]=0, & 
[{\hat h}_2,\;{\hat f}_2]=0,\\
\noalign{\medskip} [{\hat e}_1,\; {\hat f}_2]=[{\hat e}_2,\;{\hat f}_1]=0,&
{\hat e}_2^2={\hat f}_2^2=0,\\
{\hat e}_1^2{\hat e}_2-\left(q+q^{-1}\right){\hat e}_1{\hat e}_2{\hat e}_1+
{\hat e}_2{\hat e}_1^2=0,\qquad \qquad 
&{\hat f}_1^2{\hat f}_2-\left(q+q^{-1}\right){\hat f}_1{\hat f}_2{\hat f}_1+
{\hat f}_2{\hat f}_1^2=0,
\end{array}
\end{equation}        
where $\deg({\hat e}_2)=\deg({\hat f}_2)=1$ and $\deg({\hat h}_1)=\deg({\hat h}_2)=
\deg({\hat e}_1)=\deg({\hat f}_1)=0$. It is simple to note that ${\cal U}_q(sl(2))\subset 
{\cal U}_q(sl(2|1))$. The coproducts are defined by
\begin{eqnarray}
&&\Delta\left({\hat h}_1\right)={\hat h}_1\otimes 1+1\otimes {\hat h}_1,\qquad\qquad\qquad 
\Delta\left({\hat h}_2\right)={\hat h}_2\otimes 1+1\otimes {\hat h}_2,\nonumber\\
&&\Delta\left({\hat e}_1\right)={\hat e}_1\otimes q^{{\hat h}_1/2}+q^{-{\hat h}_1/2}\otimes 
{\hat e}_1,\qquad\;\;\;\Delta\left({\hat e}_2\right)={\hat e}_2\otimes q^{{\hat h}_2/2}+q^{-{\hat 
h}_2/2}\otimes {\hat e}_2, 
\nonumber\\
&&\Delta\left({\hat f}_1\right)={\hat f}_1\otimes q^{{\hat h}_1/2}+q^{-{\hat h}_1/2}\otimes 
{\hat f}_1,\qquad\; \Delta\left({\hat f}_2\right)={\hat f}_2\otimes q^{{\hat h}_2/2}+q^{-{\hat 
h}_2/2}\otimes {\hat f}_2.
\end{eqnarray}
The universal ${\cal R}$-matrix is given in ref. [27]. Note that the definition of the Hopf
superalgebra differs fromthat of the usual Hopf algebra by the supercommutativity of tenor 
product, i.e. $\left(a\otimes b\right)\left(c\otimes d\right)=(-1)^{\deg(b)\deg(c)}\left(ac
\otimes bd\right)$. For later use, we note that the fundamental representation of the algebra 
(7) is spanned by
\begin{equation}
\begin{array}{lll}
{\hat h}_1=\pmatrix{1&0&0\cr 0&-1&0\cr 0&0&0},\qquad &{\hat e}_1=\pmatrix{0&1&0\cr 0&0&0\cr 
0&0&0},\qquad & {\hat f}_1=\pmatrix{0&0&0\cr 1&0&0\cr 0&0&0},\\
{\hat h}_2=\pmatrix{0&0&0\cr 0&1&0\cr 0&0&1},\qquad & {\hat e}_2=\pmatrix{0&0&0\cr 0&0&1\cr 
0&0&0},\qquad & {\hat f}_2=\pmatrix{0&0&0\cr 0&0&0\cr 0&1&0}. 
\end{array}
\end{equation}  
  
Following the technic developped above, the ${\cal R}_{\sf h}$ matrix of the super-jordanian 
quantum superalgebra ${\cal U
}_{\sf h}(sl(2|1))$, for arbitrary representations in the two tensor product sectors, can be 
also obtainded from the ${\cal R}_q$-matrix associated with the Drinfeld-Jimbo quantum 
superalgebra ${\cal U}_q(sl(2|1))$ through a specific contraction. For simplicity and brevity, 
let us start with (fundamental irrep.) $\otimes$ (fundamental irrep.). The ${\cal R}_q$-matrix 
of ${\cal U}_q(sl(2|1))$ superalgebra in the (fund.) $\otimes$ (fund.) representation reads    
\begin{eqnarray}
R_q=\left(
\begin{array}{ccccccccc} 
q&0&0&0&0&0&0&0&0\\
\noalign{\medskip}0&1&0&q-q^{-1}&0&0&0&0&0\\
\noalign{\medskip}0&0&1&0&0&0&q-q^{-1}&0&0\\
\noalign{\medskip}0&0&0&1&0&0&0&0&0\\
\noalign{\medskip}0&0&0&0&q&0&0&0&0\\
\noalign{\medskip}0&0&0&0&0&1&0&q-q^{-1}&0\\
\noalign{\medskip}0&0&0&0&0&0&1&0&0\\
\noalign{\medskip}0&0&0&0&0&0&0&1&0\\
\noalign{\medskip}0&0&0&0&0&0&0&0&q^{-2}
\end{array}\right).
\end{eqnarray}
The non-standard ${\cal R}_{\sf h}$-matrix in the (fund.) $\otimes$ (fund.) representation is 
obtained, from (13), in the following manner:
\begin{eqnarray}
R_{\sf h}&=&\lim_{q\rightarrow 1}\left[{\sf E}^{-1}_{q}\left(\frac{{\sf h}{\hat e}_1}{q-1}
\right)_{fund.}\otimes{\sf E}^{-1}_{q}\left(\frac{{\sf h}{\hat e}_1}{q-1}\right)_{fund.}
\right]R_q\left[{\sf E}_{q}\left(\frac{{\sf h}{\hat e}_1}{q-1}\right)_{fund.}\otimes
{\sf E}_{q}\left(\frac{{\sf h}{\hat e}_1}{q-1}\right)_{fund.}\right]\nonumber\\
&=&\left(\begin{array}{ccccccccc} 
1&{\sf h}&0&-{\sf h}&{\sf h}^2&0&0&0&0\\
\noalign{\medskip}0&1&0&0&{\sf h}&0&0&0&0\\
\noalign{\medskip}0&0&1&0&0&0&0&0&0\\
\noalign{\medskip}0&0&0&1&-{\sf h}&0&0&0&0\\
\noalign{\medskip}0&0&0&0&1&0&0&0&0\\
\noalign{\medskip}0&0&0&0&0&1&0&0&0\\
\noalign{\medskip}0&0&0&0&0&0&1&0&0\\
\noalign{\medskip}0&0&0&0&0&0&0&1&0\\
\noalign{\medskip}0&0&0&0&0&0&0&0&1
\end{array}\right).
\end{eqnarray}
Similarly, using a Maple program\footnote{Our program was performed for (fund.) $\otimes$ 
(fund.), (fund.), etc.}, we obtain, for (fundamental irrep.) $\otimes$ (arbitrary 
irrep.), the follwing expression:
\begin{equation}
R_{\sf h}=\pmatrix{T& &-{\sf h}H_1+\frac{\sf h}2\left(T-T^{-1}\right)&&0\cr 
&&&&\cr
0&&T^{-1}&&0\cr
&&&&\cr
0&&0&&I},
\end{equation}
where
\begin{eqnarray}
&&T^{\pm 1}=\pm {\sf h}e_1+\sqrt{1+{\sf h}^2e_1^2},\qquad
H_1=\frac 12\left(T+T^{-1}\right)h_1=\sqrt{1+{\sf h}^2e_1^2}h_1. 
\end{eqnarray}
It is easy to verify that
\begin{eqnarray}
&&TT^{-1}=T^{-1}T=1,\qquad [H_1,T^{\pm 1}]=T^{\pm 2}-1.
\end{eqnarray} 
We note that the contraction sheme, which comprises our transformation and limiting procedure 
has furnished the ${\cal R}_h$-matrix along with a nonlinear map of the ${\sf h}$-Borel 
subalgebra on its classical counterpart. Following ref. \cite{P21}, let us introduce the 
generator
\begin{equation}
F_1=f_1-\frac{{\sf h}^2}4e_1\left(h_1^2-1\right).
\end{equation} 
We then show that
\begin{eqnarray}    
&&[T^{\pm 1},F_1]=-\frac{\sf h}2\biggl(H_1T^{\pm 1}+T^{\pm 1}H_1\biggr),\qquad 
[H_1,F_1]=-{1\over 2}\biggl(TF_1+F_1T+T^{-1}F_1+F_1T^{-1}\biggr).
\end{eqnarray}
The coproducts $\Delta$, the counit $\epsilon$ and the antipode $S$ of $\{H,T,T^{-1},Y\}$ read 
\cite{P4}  
\begin{eqnarray}
&&\Delta(H_1)=H_1\otimes T+T^{-1}\otimes H_1,\qquad\qquad \Delta(T^{\pm 1})=T^{\pm 1}\otimes 
T^{\pm 1},\qquad\qquad  \Delta(F_1)=F_1\otimes T+T^{-1}\otimes F_1,\nonumber\\
&&S(H_1)=-TH_1T^{-1},\qquad S(T^{\pm 1})=T^{\mp 1},\qquad S(F_1)=-TF_1T^{-1},\nonumber\\
&&\epsilon(H_1)=\epsilon(F_1)=0,\qquad\epsilon(T^{\pm 1})=1.
\end{eqnarray}
This implies that the Ohn's structure follows from the bosonic generators $\{h_1,e_1,f_1\}$. 
The algebraic properties (99) exhibits clearly the embedding of ${\cal U}_{\sf h}(sl(2))$ in 
${\cal U}_{\sf h}(sl(2|1))$. At this level it is natural to ask the following question. {\it 
How to complete the ${\cal U}_{\sf h}(sl(2|1))$ superalgebra?}.    
To complete now the ${\cal U}_{\sf h}(sl(2|1))$ superalgebra, we introduce the following 
${\sf h}$-deformed fermionic root generator:
\begin{eqnarray}
&&H_2=h_2-\frac{{\sf h}^2}2e_1^2h_1,\qquad 
E_2=e_2-\frac{{\sf h}^2}4e_1e_3\left(2h_1+1\right),\qquad
F_2=f_2,\nonumber\\
&&H_3=h_3+\frac{{\sf h}^2}2e_1^2h_1,\qquad
E_3=e_3,\qquad F_3=f_3+\frac{{\sf h}^2}4e_1f_2\left(2h_1+1\right)
\end{eqnarray}   
Let us note that there exist a $\CC$-algebra automorphism of ${\cal U}_{\sf h}(sl(2|1))$ 
such that
\begin{equation}
\begin{array}{lll}
\Phi\left(T^{\pm 1}\right)=T^{\pm 1},\qquad\qquad &\Phi\left(F_1\right)=F_1,\qquad \qquad 
&\Phi\left(H_1\right)=H_1, \\
\Phi\left(E_2\right)=F_3,\qquad\qquad &\Phi\left(F_2\right)=-E_3,\qquad\qquad 
&\Phi\left(H_2\right)=-H_3, \\
\Phi\left(E_3\right)=-f_2,\qquad\qquad &\Phi\left(F_3\right)=E_2,\qquad\qquad 
&\Phi\left(H_3\right)=-H_2.
\end{array}
\end{equation}
(For ${\sf h}=0$, this automorphism reduces to the classical one described by (89). The 
expressions (96), (97), (98) and (101) define a realisation of the 
super-jordanian subalgebra ${\cal U}_{\sf h}(sl(2|1))$ with the classical generators via a 
nonlinear map. Other invertible maps relating the super-jordanian and the classical generators 
may also be considered. Our construction leads to the following results (We have quoted here 
only the final results): The nonstandard (super-jordanian) enveloping superalgebra ${\cal 
U}_{\sf h}(sl(2|1))$ is an associative superalgebra over $\CC$ generated by $\{H_1,\;T,\;T^{-1}
,\;F1,\;H_2,\;E_2,\;F_2,\;H_3,\;E_3,\;F3\}$, along with (99), the commutation relations [24]
\begin{eqnarray}
&&[H_1,\;H_2]=-\frac 14\left(T-T^{-1}\right)H_1,\qquad
[H_1,\;H_3]=\frac 14\left(T-T^{-1}\right)H_1,\qquad [H_2,\;H_3]=0,\nonumber\\
&& [H_1,\;E_2]=-\frac 12E_2\left(T+T^{-1}\right)-\frac{\sf h}4\left(H_1\left(T-T^{-1}\right)+
\left(T-T^{-1}\right)H_1\right)E_3,\nonumber\\
&&[H_1,\;F_2]=\frac 12\left(T+T^{-1}\right)F_2,\qquad\qquad [H_1,\;E_3]=\frac 12\left(T+T^{-1}
\right)E_3,\nonumber\\ 
&&[H_1,\;F_3]=\frac 12F_3\left(T+T^{-1}\right)+\frac{\sf h}4\left(H_1\left(T-T^{-1}\right)+
\left(T-T^{-1}\right)H_1\right)F_2,\nonumber\\ 
&&[H_2,\;T]=-\frac 14\left(T^3-T^{-1}\right),\qquad\qquad 
[H_2,\;T^{-1}]=-\frac 14\left(T^{-3}-T\right),\nonumber\\
&&[H_2,\;F_1]=\frac 14\left(T+T^{-1}\right)^2F_1-\frac{\sf h}4\left(T-T^{-1}\right)H_1^2
-\frac{\sf h}4\left(T^2-T^{-2}\right)H_1-\frac{\sf h}{16}\left(T^2-T^{-2}\right)\left(T+
T^{-1}\right),\nonumber\\
&& [H_2,\;E_2]=\frac{\sf h}{16}\left(T-T^{-1}\right)\left(T^2-T^{-2}\right)E_3
+\frac 18\left(T-T^{-1}\right)^2E_2,\nonumber\\
&& [H_2,\;F_2]=-\frac 18\left(T-T^{-1}\right)^2F_2,\qquad 
[H_2,\;E_3]=-\frac 18\left(T^2+6+T^{-2}\right)E_3,\nonumber\\
&& [H_2,\;F_3]=\frac 18\left(T^2+6+T^{-2}\right)F_3-\frac{\sf h}{16}\left(T^2-T^{-2}\right)
\left(T+T^{-1}\right)F_2,\nonumber\\ 
&&[H_3,\;T]=\frac 14\left(T^3-T^{-1}\right),\qquad\qquad 
[H_3,\;T^{-1}]=\frac 14\left(T^{-3}-T\right),\nonumber\\
&&[H_3,\;F_1]=\frac 14\left(T+T^{-1}\right)^2F_1+\frac{\sf h}4\left(T-T^{-1}\right)H_1^2
+\frac{\sf h}4\left(T^2-T^{-2}\right)H_1+\frac{\sf h}{16}\left(T^2-T^{-2}\right)\left(T+
T^{-1}\right),\nonumber\\
&& [H_3,\;E_2]=-\frac 18\left(T^2+6+T^{-2}\right)E_2-\frac{\sf h}{16}\left(T^2-T^{-2}\right)
\left(T+T^{-1}\right)E_3,\nonumber\\ 
&&[H_3,\;F_2]=\frac 18\left(T^2+6+T^{-2}\right)F_2,\qquad
[H_3,\;E_3]=\frac 18\left(T-T^{-1}\right)^2E_3,\nonumber\\
&&[H_3,\;F_3]=\frac{\sf h}{16}\left(T-T^{-1}\right)\left(T^2-T^{-2}\right)F_2
-\frac 18\left(T-T^{-1}\right)^2F_3,\nonumber\\
&&[E_2,\;F_2]=H_2-\frac 1{16}\left(T-T^{-1}\right)^2-\frac{\sf h}4\left(T-T^{-1}\right)E_3F_2,\nonumber\\
&&[E_3,\;F_3]=H_3+\frac 1{16}\left(T-T^{-1}\right)^2+\frac{\sf h}4\left(T-T^{-1}\right)F_2E_3,
\qquad [T^{\pm 1},\;F_2]=0,\nonumber\\
&&[E_2,\;F_1]=\frac{\sf h}4\left(T-T^{-1}\right)E_2+\frac{\sf h}2\left(T-T^{-1}\right)E_3F_1
-\frac{{\sf h}^2}4E_3H_1^2\nonumber\\
&&\phantom{[E_2,\;F_1]=}-\frac{3{\sf h}^2}8\left(T+T^{-1}\right)E_3H_1-
\frac{{\sf h}^2}2E_3-\frac{15{\sf h}^2}{64}\left(T-T^{-1}\right)^2E_3,\nonumber\\
&&E_2^2=\frac{\sf h}4\left(T-T^{-1}\right)E_3E_2,\qquad
F_2^2=E_3^2=0,\qquad 
F_3^2=-\frac{\sf h}4\left(T-T^{-1}\right)F_2F_3,\nonumber\\
&&[T^{\pm 1},\;E_2]=\pm \frac{\sf h}2\left(T^{\pm 2}+1\right)E_3,\qquad 
[F_2,\;F_1]=F_3,\qquad [T^{\pm 1},\;E_3]=0,\nonumber \\
&&[F_3,\;F_1]=\frac{\sf h}4\left(T-T^{-1}\right)F_3-\frac{\sf h}2\left(T-T^{-1}\right)F_2F_1
+\frac{{\sf h}^2}4F_2H_1^2\nonumber\\
&&\phantom{[E_2,\;F_1]=}+\frac{3{\sf h}^2}8\left(T+T^{-1}\right)F_2H_1+
\frac{{\sf h}^2}2F_2+\frac{15{\sf h}^2}{64}\left(T-T^{-1}\right)^2F_2,\nonumber\\
&&[E_2,\;E_3]=0,\qquad\qquad [F_2,\;F_3]=0,\qquad [F_1,\;E_3]=E_2,\nonumber\\
&&[T,\;F_3]=-\frac{\sf h}2\left(T^2+1\right)F_2,\qquad 
[T^{-1},\;F_3]=\frac{\sf h}2\left(T^{-2}+1\right)F_2,\qquad 
[F_2,\;E_3]=\frac1{2{\sf h}}\left(T-T^{-1}\right),\nonumber\\
&&[F_3,\;E_2]=F_1-\frac{\sf h}4\left(T-T^{-1}\right)F_2E_2+\frac{\sf h}4\left(T-T^{-1}\right)
E_3F_3-\frac{\sf h}8\left(T-T^{-1}\right)H_1^2-\frac{\sf h}8\left(T^2-T^{-2}\right)H_1
\nonumber\\
&&\phantom{[F_3,\;E_2]=}-\frac{\sf h}{16}H_1\left(T^2-T^{-2}\right)-
\frac{7{\sf h}}{128}\left(T-T^{-1}\right)^3. 
\end{eqnarray}
The $\ZZ_{2}$-grading in ${\cal U}_{\sf h}(sl(2|1))$ is uniquely defined by the requirement 
that the only odd generators are $E_2$, $F_2$, $E_3$ and $F_3$, i.e. $\deg\;(H_1)=\deg\;(H_2)
=\deg\;(H_3)=\deg\;(T)=\deg\;(T^{-1})=\deg\;(F_1)=0$ and $\deg\;(E_2)=\deg(F_2)=\deg\;(E_3)=
\deg(F_3)=1$. It is obvious that as ${\sf h}\rightarrow 0$, we have $\left(E_2,F_2,H_2,E_3,F_3,
H_3\right)\rightarrow\left(e_2,f_2,h_2,e_3,f_3,h_3\right)$. 
 
We turn now to the coalgebraic structure. The nonstandard (super-jordanian) quantum 
enveloping superalgebra ${\cal U}_{\sf h}(sl(2|1))$ admits a Hopf structure with coproducts, 
antipodes and counits determined by (100) and [24]
\begin{eqnarray}
&&\Delta\left(E_2\right)=E_2\otimes T^{1/2}+T^{-1/2}\otimes E_2+\frac{\sf h}4T^{-1}E_3\otimes 
\left(T^{-1/2}H_1+H_1T^{-1/2}\right)-\frac{\sf h}4\left(T^{1/2}H_1+H_1T^{1/2}\right)\otimes 
TE_3,\nonumber\\
&&\Delta\left(F_2\right)=F_2\otimes T^{-1/2}+T^{1/2}\otimes F_2,\qquad 
\Delta\left(E_3\right)=E_3\otimes T^{-1/2}+T^{1/2}\otimes E_3,\nonumber\\
&&\Delta\left(F_3\right)=F_3\otimes T^{1/2}+T^{-1/2}\otimes F_3-\frac{\sf h}4T^{-1}F_2\otimes 
\left(T^{-1/2}H_1+H_1T^{-1/2}\right)+\frac{\sf h}4\left(T^{1/2}H_1+H_1T^{1/2}\right)\otimes 
TF_2,\nonumber\\
&&\Delta\left(H_2\right)=H_2\otimes 1+1\otimes H_2+\frac 14TH_1\otimes\left(1-T^2\right)+
\frac 14\left(1-T^{-2}\right)\otimes T^{-1}H_1,\nonumber\\
&&\Delta\left(H_3\right)=H_2\otimes 1+1\otimes H_2-\frac 14TH_1\otimes\left(1-T^2\right)-
\frac 14\left(1-T^{-2}\right)\otimes T^{-1}H_1,\nonumber\\
&&S\left(E_2\right)=-E_2-\frac{\sf h}2\left(T+T^{-1}\right)E_3,\qquad 
S\left(F_2\right)=-F_2,\qquad S\left(E_3\right)=-E_3,\nonumber\\
&&S\left(F_3\right)=-F_3+\frac{\sf h}2\left(T+T^{-1}\right)E_3,\qquad 
S\left(H_2\right)=-H_2+\frac 12\left(T^{-2}-1\right),\qquad 
S\left(H_3\right)=-H_3-\frac 12\left(T^{-2}-1\right),\nonumber\\
&&\epsilon\left(H_2\right)=\epsilon\left(H_3\right)=\epsilon\left(E_2\right)=\epsilon
\left(F_2\right)=\epsilon\left(E_3\right)=\epsilon\left(F_3\right)=0.  
\end{eqnarray}
All the Hopf superalgebra axioms can be verified by direct calculations. The universal 
${\cal R}_{\sf h}$-matrix of ${\cal U}_{\sf h}(sl(2|1))$ has the following form [24]:  
\begin{eqnarray}
&& {\cal R}_{\sf h}=\exp\biggl(-{\sf h}X_1\otimes TH_1\biggr)\exp\biggl({\sf h}TH_1\otimes 
X_1\biggr),  
\end{eqnarray}
where $X_1={\sf h}^{-1}\ln T$. The element (105) coincides with the pure 
${\cal U}_{\sf h}(sl(2))$ universal ${\cal R}_{\sf h}$-matrix [26].

\section{Conclusion}

In general, a class of nonlinear invertible maps exists relating the Jordanian quantum 
algebras and their classical analogues. Here we have used a particular maps realizing 
Jordanian ${\cal U}_{\sf h}(sl(2))$, ${\cal U}_{\sf h}(sl(3))$, ${\cal U}_{\sf h}(sl(N))$,
${\cal U}_{\sf h}(osp(2|1))$ and ${\cal U}_{\sf h}(sl(2|1))$. As a result a result of choice 
of the basis, via the maps described earlier, the algebraic commutations relations are 
deformed. On benifict of our procedure is that our expressions of the coalgebraic structure 
are considerably simpler than those found elswhere [8-20].           

\medskip
\medskip
\large{Comment}: Talk given by MB. ZAHAF to the Seventh Constantine High 
Energy Physics School 
(Theoretical Physics ans Cosmology), 3-7 April 2004, Constantine (Algeria).

\end{document}